\input amstex
\documentstyle{amsppt}
\magnification1200 
\tolerance=1000
\def\n#1{\Bbb #1}

\def\fr{\hbox{fr }}

\def\lim{\hbox{lim }}
\def\im{\hbox{im }}
\def\id{\hbox{id }}
\def\ad{\hbox{ad }}

\def\opr{\hbox{def}}
\def\Res{\hbox{Res }}
\def\dim{\hbox{dim }}

\def\rank{\hbox{ rank }}
\def\Corr{\hbox{ Corr }}
\def\diag{\hbox{ diag }}

\def\suchthat{\hbox{ such that }}

\topmatter
\title 
Relations between conjectural eigenvalues of Hecke operators on submotives of Siegel varieties
\endtitle 
\author
Dmitry Logachev 
\endauthor
\NoRunningHeads 
\endtopmatter
\document 

{\bf Abstract. }
\medskip
There exist conjectural formulas on relations between $L$-functions of submotives of Shimura varieties and automorphic representations of the corresponding reductive groups, due to Langlands --- Arthur. In the present paper these formulas are used in order to get explicit relations between eigenvalues of $p$-Hecke operators (generators of the $p$-Hecke algebra of $X$) on cohomology spaces of some of these submotives, for the case $X$ is a Siegel variety. Hence, this result is conjectural as well: methods related to counting points on reductions of $X$ using the Selberg trace formula are not used. 

It turns out that the above relations are linear, their coefficients are polynomials in $p$ which satisfy a simple recurrence formula. The same result can be easily obtained for any Shimura variety. 

This result is an intermediate step for a generalization of the Kolyvagin's theorem of finiteness of Tate -- Shafarevich group of elliptic curves of analytic rank 0 , 1 over $\n Q$, to the case of submotives of other Shimura varieties, particularly of Siegel varieties of genus 3. 

Idea of the proof: from one side, the above formulas of Langlands --- Arthur give us (conjectural) relations between Weil numbers of a submotive. From another side, the Satake map permits us to transform these relations between Weil numbers into relations between eigenvalues of $p$-Hecke operators on $X$. 

The paper contains also a survey of some related questions, for example explicit finding of the Hecke polynomial for $X$, and (Appendix) tables for the cases $g=2,3$. 
\medskip
Version of Feb. 19, 2002

\newpage
{\bf Introduction.}\nopagebreak
\medskip
The purpose of the present paper is to show that starting from some standard conjectures of Langlands --- Arthur, a chain of elementary calculations leads to a simply-formulated and non-expected result on  relations between eigenvalues of $p$-Hecke operators on a Shimura variety (which  hence depends on these conjectures). 

Namely, let $\Pi$ be a stable global packet of automorphic representations of a reductive group $G$ ($G$ corresponds to a Shimura variety $X$). Attached to $\Pi$ is a parabolic subgroup $P$ of $G$. Let $\Cal M$ be an irreducible constituent of a submotive conjecturally attached to $\Pi$ (see (0.1) below). We denote by $E$ the field of coefficients of $\Cal M$. $H^i(\Cal M)$ is a module over the $p$-Hecke algebra $\n H(G)=\n H_p(G)$ of $X$, $p$ is a fixed prime. Any $t \in \n H(G)$ acts on $H^i(\Cal M)$ by multiplication by an element $\goth m_{\Cal M}(t)$ of $E$. If $\{t_i\}$ is a set of generators of $\n H(G)$ then numbers $\goth m_{\Cal M}(t_i)$ can satisfy some relations. 

The main result of the present paper is the Theorem 4.4. We fix the type of a Shimura variety (the level is arbitrary) and a parabolic subgroup $P$. We find the set of relations between $\goth m_{\Cal M}(t_i)$ (depending on $P$ only) such that if the above conjectures are true then these relations are satisfied for all submotives $\Cal M$ corresponding to all stable $\Pi$ corresponding to $P$. 

Really, the first steps of the calculations of the present paper are made only for Siegel varieties, and the last step only for submotives corresponding to one type of $P$ (the simplest no-trivial). This restriction is not of principle: a reader can easily get analogous results for any Shimura variety and any type of submotives, including the ones that correspond to non-stable packets. The author is interested in the case $g=3$, because the main result of the present paper can be applied for a generalization of the Kolyvagin's theorem to the case of submotives of Siegel varieties of genus 3.  

The paper contains also some calculations that are not logically necessary for the proof of the main theorem, but can be used for references or for the further development of the subject (Sections 2.2 - 2.4, 3.4 and small parts of other sections). Some tables for genii 2 and 3 are given in the Appendix, which is written by the same reason. 
\medskip
{\bf Acknowledgments.} The author is grateful to A. Andrianov and M. Borovoi for important advice on the subject of this paper. 
\medskip
\medskip
{\bf 0. Idea of the proof.}\nopagebreak
\medskip
Recall that in order to define a Shimura variety $X$ of a fixed level, we must fix a reductive group $G$ over $\n Q$, a map $h: \Res_{\n C/\n R}(G_m)  \to G$ over $\n R$, and a level subgroup $\Cal K \subset G(\n A_f(\n Q))$; these data must satisfy some conditions ([D]). Throughout the paper we consider only the case when $p$ does not divide the level, i.e. $\Cal K \supset G(\n Z_p)$. We shall consider only the case of Siegel varieties, i.e. from here we let $G=GSp_{2g}$ over $\n Q$. Further, we must choose a compactification of $X$ and a type of cohomology. 

Really, all subsequent considerations depend only on $G$ and $h$ and hence does not depend on level, compactification and the type of cohomology. 

We fix a Borel subgroup $B$ of $G$ and consider all intermediate parabolic subgroups $P$ of $G$, $B\subset P \subset G$. There is a 1 - 1 correspondence between the set of archimedean cohomological A-parameters of $G$ and the set of such $P$ ([AJ]; [A]; [BR], Section 4.1). We denote by $\Pi_P$ the packet of automorphic representations of $G(\n R)$ corresponding to the archimedean cohomological A-parameter corresponding to $P$ ([BR], Section 4.2). 

Attached to $P$ and $X$ is a set (indexed by $k$) of stable global packets of automorphic representations of $G(\n A_{\n Q})$. This set clearly depends on the level of $X$. We denote the $k$-th packet by $\Pi^{glob}_P (k)$. Let $\pi \in \Pi^{glob}_P (k)$ be a representation, $\pi=\pi_{\infty}\otimes\pi_f$ its decomposition on archimedean and finite part, and $\pi_f = \otimes_{l} \pi_{l}$ ($l =$ prime) the decomposition of $\pi_f$. We consider only such $\pi$ for which $\pi_p$ is non-ramified (see Step 3 below for a description of $\pi_p$). We have $\pi_{\infty} \in \Pi_P$, and for any other $\pi'_{\infty}\in \Pi_P$ the representation $\pi' \overset{\opr}\to{=}\pi'_{\infty}\otimes\pi_f$ also belongs to $\Pi^{glob}_P (k)$. 

Conjecturally, $\forall k$ there exists a submotive $\Cal M_P(k)$ (reducible unless $P=B$) such that $$L_p(\pi, r,s)=L_p(\Cal M_P(k),s)\eqno{(0.1)}$$ where $r: ^L G\to GL(W)$ is a finite-dimensional representation defined in [BR], 5.1; Weil numbers of $\Cal M_P(k)$ and hence its $L$-function are considered with respect to $E$, and $L_p(\pi, r,s)=L(\pi_p, r,s)$ is the local $L$-function (particularly, it does not change if we change $\pi$ by $\pi'$). 

From now we fix some $k$ and we denote $\Cal M_P(k)$ simply by $\Cal M_P$. It is the main object of the present paper. Since it is conjectural, all subsequent theorems should be understand as follows: 
\medskip
Let $\Cal M_P$ be a motive such that (0.1) is satisfied. Then ... (statement of the theorem). 
\medskip
Some references for properties of $\Cal M_P$: [BR], Sections 4.3, 5, 7; [A], Section 9, and the present paper, section 3. Most of these properties are not necessary for the proof of the main theorem. Here we mention only that (0.1) implies that $\sum_i h^i(\Cal M_P)=\dim r$; for $G=GSp_{2g}$ this number is $2^g$. See also Remark 0.1 below. 
\medskip
Now let $\Cal M \subset \Cal M_P$, $\goth m_{\Cal M}$, $E$ be as in the Introduction. Multiplication by elements of $E$ gives us an inclusion $i_{\Cal M}: E \to \n C$. It is clear that the composite map $i_{\Cal M}\circ\goth m_{\Cal M}:\n H(G) \to \n C$ does not depend on ${\Cal M}$ but only on ${\Cal M}_P$. So, we denote the map $i_{\Cal M}\circ\goth m_{\Cal M}$ by $\goth m=\goth m_P$. The main result of the present paper is the finding of relations between numbers $\goth m(t_i)$. 
\medskip
Step 1. We use the following notation: if $b_1,\dots, b_g$ is a set of numbers and $I$ is a subset of the set $\{1,\dots,g\}$ then we denote $b_I=\prod_{i\in I}b_i$. In (2.7) we recall a (well-known; see, for example, [BR], Sect. 5.1, Example 3) proof of a fact that there exists a set of numbers $a_0, b_1, b_2,\dots, b_g\in \n C$ such that eigenvalues of $r(\theta_{\pi_p})$ (where $\theta_{\pi_p}$ is a Langlands element of $\pi_p$) have the form $a_0b_I$ where $I$ runs over the set of all $2^g$ subsets of $(1,\dots,g)$. Hence, (0.1) means that the Weil numbers of $\Cal M_P$ have the same form. 
\medskip
Step 2. We shall show in Sect. 4 (as a result of calculations of Sect. 3) that for a fixed $P$ numbers $a_0, b_1, b_2,\dots, b_g$ satisfy some relations depending only on $P$ (Prop. 4.3). 
\medskip
Step 3. Finally, using explicit formulas of Satake map (Sect. 2), we shall show that relations of Prop. 4.3 give us relations between numbers $\goth m(t_i)$. 
\medskip
Now we describe steps 2 and 3 in more details. 
\medskip
{\bf Step 2a.} Firstly, we recall the description of the set of parabolic subgroups $P$ under consideration: there are 2 types of such subgroups, and subgroups of each type are parametrized by the set of ordered partitions of $g$, i.e. the set of representations of $g$ as a sum 
$$g=\goth b_1+\goth b_2+\dots+\goth b_k\eqno{(0.2)}$$ 
($\goth b_i\ge 1$, the order is essential), or, the same, the set of sequences 
$$0=m_1<m_2<\dots<m_{k+1}=g$$
where $m_i=\goth b_1+\cdots \goth b_{i-1}$.
\medskip
{\bf Step 2b.} Using formulas of [BR], Sect. 4, we describe explicitly in Sect. 3.3 the set of archimedean cohomological representations belonging to $\Pi_P$. Namely, let $P$ be of type 1 given by (0.2). We denote by $\goth C$ the set of sequences $\goth c=(c_1,\dots,c_k)$ such that $\forall j=1,\dots, k$ $0\le c_j\le \goth b_j$. We have: $\Pi_P$ is isomorphic to $\goth C$ factorized by the equivalence relation $(c_1,\dots,c_k) \sim (\goth b_1-c_1,\dots,\goth b_k-c_k)$. For $P$ of type 2 the result is the same, but $c_1$ is omitted. The representation corresponding to $\goth c \in \goth C$ is denoted by $\pi_{\goth c}$. 
\medskip
{\bf Step 2c.} Now we use formulas of [BR], Sect. 4 for the dimensions of $H^{i,j}(\goth g,K_c;\pi_{\goth c})$, $\goth{g=gsp}_{2g}$. We consider for all $i=1,\dots,k$ the set $\goth S(c_i,\goth b_i)$ of all subsets of order $c_i$ of the set $\{1,\dots,\goth b_i\}$, and we denote 
$$\goth S(\goth c, P)=\prod_{i=1}^k \goth S(c_i,\goth b_i)\hbox{ or }\goth S(\goth c, P)=(\n Z/2\n Z)^{\goth b_1}\prod_{i=2}^k \goth S(c_i,\goth b_i)\eqno{(0.3)}$$ for $P$ of type 1 and 2 respectively (this is the set of representatives of minimal length for the cosets $\Omega(T,M/\Omega(T,M\cap wKw^{-1})$ in notations of [BR], Sect. 4.3). For $\rho\in \goth S(\goth c, P)$ an explicit formula for the length $l(\rho)$ is given in (3.5.1). 

Further, (3.5.1) give us relations between Weil numbers of $\Cal M_P$. The exact formula for these relations is given in (4.1). Really, (4.1) is a corollary of a stronger proposition 4.3. 
\medskip
**********************************************************
\medskip
The following Remark 0.1 is not a part of the text of the paper. It is included in order 

(1) To help to a reader to understand the notations; 

(2) It contains some questions which I cannot answer, for example: (a) is condition (0.1.2) sufficient? (b) Let us consider a decomposition of the right hand side of (0.1.4) coming from the theorem on structure of finite-dimensional $\goth{sl}_2$-modules. Does this decomposition come from a decomposition of $\Cal M_{P,i}$ itself? 

I would be very grateful to a reader for answers. 

Since the contents of the Remark 0.1 are not used in the present paper, possible errors have not influence on the proof of the main theorem. 
\medskip
{\bf Remark 0.1.} A sketch of the description of the structure of $\Cal M_P$ ([A], [BR]).
\medskip
The space generated by $\goth S(\goth c, P)$ is isomorphic to $\bigoplus_{i,j}H^{i,j}(\goth g,K_c;\pi_{\goth c})$, so it is a $\goth{sl}_2$-module with a Hodge structure. We denote this $\goth{sl}_2$-module by $Lie(\goth S(\goth c, P))$. There are 2 numbers $p_{\goth c}, q_{\goth c}$ associated to each $\goth c$ (see [BR], Section 4.3 for a formula for them, and (3.4) for explicit values). 2 basis elements of minimal weight in $\goth S(\goth c, P)$ have Hodge numbers $h^{p_{\goth c}, q_{\goth c}}=h^{q_{\goth c}, p_{\goth c}}=1$ (case $G=GSp$, $p_{\goth c}\ne q_{\goth c}$; formulas for other Hodge numbers are given for example in [A], Section 9). 

There exists a partition of $\goth C$: $$\goth C=\bigcup_{i\in \goth I}\goth C_{i}\eqno{(0.1.1)}$$ (the union is disjoint) which gives rise to a descomposition of $\Cal M_P$ as a direct sum of submotives. For $\goth c_1$, $\goth c_2$ a necessary condition to belong to one $\goth C_{i}$ is the following: 

{\bf(0.1.2)} For $j=1,2$ numbers $p_{\goth c_j}+q_{\goth c_j}$ coincide and $\goth{sl}_2$-modules $Lie(\goth S(\goth c_j, P))$ (but not their Hodge structures!) are isomorphic.

Attached to (0.1.1) is a motive decomposition $$\Cal M_P=\bigoplus_{i\in \goth I}\Cal M_{P,i}\eqno{(0.1.3)}$$ having the following property: $H^*(\Cal M_{P,i})$ has a natural structure of $\goth{sl}_2$-module, and we have an isomorphism of $\goth{sl}_2$-modules with Hodge structures: 
$$H^*(\Cal M_{P,i})=\bigoplus_{\goth c\in \goth C_i} Lie(\goth S(\goth c, P))\eqno{(0.1.4)}$$
and analogously for their components of any fixed weight. It is known that a descomposition of $\Cal M_{P,i}$ in a direct sum indexed by $\goth c\in \goth C_i$ --- like in (0.1.4) --- does not exist. Clearly (0.1.3), (0.1.4) give us a description of Hodge numbers of $\Cal M_P$ and primitive elements in its cohomology groups. 

See also Appendix, 8 for some explicit properties of $\Cal M_P$, where $P$ is of two simplest types.
\medskip
**********************************************************
\medskip
{\bf Step 2d.} To complete Step 2, we must use results of steps 1 and 2c in order to find relations between numbers $a_0, b_1, b_2,\dots, b_g$. These relations are the following (Proposition 4.3): 
\medskip
{\bf(0.4)} $P$ of the first type: $b_{m_i+1}$ are free variables, $b_{m_i+c}=p^{c-1}b_{m_i+1}$ ($c=1,\dots,\goth b_i$), and $a_0$ is defined by the equality $a_0^2\prod b_i=p^{g(g+1)/2}$. 

$P$ of the second type: $b_i=p^i$ for $i=1,\dots, \goth b_1$, $b_{m_i+1}$ ($2\le i\le k$) are free variables, $b_{m_i+c}$ and $a_0$ are like above. 
\medskip
\medskip
{\bf Step 3.} The $p$-Hecke algebra $\n H(G)$ is the ring of polynomials whose generators are denoted by $\tau_{p,*}$, $*=\emptyset, 1, \dots,g$:   
$\n H(G) = \n Z[\tau_p, \tau_{p,1}, \dots , \tau_{p,g}]$. Let $\chi : T(\n Q_p) \to \n C^*$ be a nonramified character such that $\pi_{\chi}=\pi_p$ where $\pi_p$ is the $p$-part of $\pi$ and 
$\pi_{\chi}: G(\n Q_p) \to GL(\Cal V)$ is the parabolically 
induced representation. $\chi$ does not depend on $\pi \in \Pi^{glob}_P (k)$. $\pi_{\chi}$ defines an action of $\n H(G)$ on a 1-dimensional subspace of $\Cal V^{G(\n Z_p)}$ and hence a homomorphism $\alpha_G(\chi): \n H(G) \to \n C$. Obviously $\alpha_G(\chi)=\goth m$, hence in order to find relations between numbers $\goth m(\tau_{p,*})$ we need to represent $\alpha_G(\chi)(\tau_{p,*})$ as polynomials in $a_0, b_1, b_2,\dots, b_g$ and to use (0.4). 
\medskip
To solve this problem we use 
\medskip
(a) the Satake map $S:\n H(G) \to \n H(T)$ where $\n H(T)\overset{i}\to{\hookrightarrow} \n Z[U_i^{\pm 1}, V_i^{\pm 1}]$ is the Hecke algebra of a maximal torus $T$ of $G$; 
\medskip
(b) an explicit expression for a Langlands element $\theta(\pi_{\chi})\in \hat T \subset \; ^L G$ given in (2.7.1), and a decomposition of $r|_{\hat T}$ as a sum of characters of $\hat T$ (Section 2.6). 
\medskip
Explicit formulas for $S(\tau_{p,*})$ are given in Section 1 ((1.2.1), (1.3.1), (1.5.1)). Further, there exists a map $\alpha_T(\chi): \n Z[U_i^{\pm 1}, V_i^{\pm 1}] \to \n C$ such that $\alpha_G(\chi)=\alpha_T(\chi)\circ i\circ S$. (0.1), (2.5.1) and (2.7.1) show us that $\alpha_T(\chi)(V_i)=a_0^{1/g}$, $\alpha_T(\chi)(U_i)=a_0^{1/g}b_i$.

Using explicit formulas for $i\circ S$ (Section 1), we can represent $\alpha_G(\chi)(\tau_{p,*})$ as polynomials in $\alpha_T(\chi)(U_i)$, $\alpha_T(\chi)(V_i)$, i.e. as polynomials in $a_0, b_1, b_2,\dots, b_g$ (2.7.4). The final result follows immediately from (0.4) and (2.7.4). 
\medskip
{\bf Structure of the paper.} In Section 1.1 we recall the definition of Satake maps $S_G$, $S_T$ and define generators of Hecke algebras $\n H(G)$, $\n H(M_s)$. In 1.2, 1.3 we find explicitly $S_G$ of these generators. Remark 1.4 is used only for a proof that 2 methods of finding of Hecke polynomial give the same result. Remark 1.5 gives a slightly other method of description of Satake map; some notations of 1.5 will be used later. 

Section 2.1 contains a definition of the induced representation and of the corresponding map $\alpha_G(\chi): \n H(G) \to \n C$. Sections 2.2 - 2.4 are of survey nature: they contain explicit formulas for $\alpha_G(\chi)$ using counting of cosets. A formula for $\alpha_G(\chi)$ that will be really used in future is given in 2.5. In 2.6 we recall properties of the map $r$ which is used to define the $L$-function of $\Cal M$, and in 2.7 we get an expression for Weil numbers of $\Cal M$. 

In 3.1 we recall the definition of parabolic subgroups of $G$ and related groups. Contents of other subsections 3.2 --- 3.5 correspond to their titles. Finally, Section 4 contains the end of the proof. 
\medskip
\medskip
{\bf 1. Explicit description of Satake map. }
\medskip
{\bf 1.1.} References: [AZh], [FCh]. We let: $T\subset G$ is a torus of diagonal matrices; $M_s=\left\{\left( \matrix A & 0 \\ 0 & (A^t)^{-1} \endmatrix \right)\right\}\subset G$.

Here we consider elements of $p$-Hecke algebras $\n H(\goth G)$ ($\goth G=G,M_s,T$) as linear combinations of double cosets of $\goth G(\n Z_p)$. There are inclusions $\n H(G) \subset \n H(M_s) \subset \n H(T)$ defined by Satake maps denoted by $S_G$, $S_T$ respectively (see [S], [FCh]).

We need the following matrices: 

$T_p = \left( \matrix 1 & 0 \\ 0 & p \endmatrix \right) $, 
entries are $g \times g$-matrices;

$T_{p,i} = \left( \matrix 1 & 0 &   0 & 0 
                \\ 0 & p &   0 & 0
                \\ 0 & 0 & p^2 & 0 
                \\ 0 & 0 &   0 & p \endmatrix \right) $, 
diagonal entries are $g - i \times g - i$, $i \times i$, $g - i
\times g - i$, $i \times i$-matrices,  $ i = 0,\dots,g$. 

We denote the double cosets $G(\n Z_p) T_p G(\n Z_p)$, $G(\n Z_p) T_{p,i} G(\n Z_p)$ 
(= elements of 
$\n H(G)$) by $\tau_p$, $\tau _{p,i}$ respectively. 
It is known that $\n H(G)$ is the ring of polynomials: 
$\n H(G) = \n Z[\tau_p, \tau_{p,1}, \dots , \tau_{p,g}]$. 
Now we need matrices 
$$F_{p,i} = F_i = \left( \matrix 1 & 0 & 0 & 0
                            \\ 0 & p & 0 & 0
                            \\ 0 & 0 & p & 0
                            \\ 0 & 0 & 0 & 1 \endmatrix 
							\right) $$
where diagonal entries are $g - i \times g - i$, $i \times i$, $g - i
\times g - i$, $i \times i$-matrices,  $ i = 0,\dots,g$. 

We denote the corresponding elements $M_s(\n Z_p)F_iM_s(\n Z_p)$ of $\n H(M_s)$ by $\Phi_i$. 
\medskip
Let us recall the definition of the Satake map $S_G$. Here we consider for $*=G$ or $M_s$ an element $f \in \n H(*)$ as a $*(\n Z_p)$-biinvariant function on $*(\n Q_p)$; a function associated to a double coset is its characteristic function. $S_G(f)$ is defined completely by its values on elements 
$X \in M_s(\n Q_p)$ of the form $X = \diag (p^{a_1}, 
\dots , p^{a_g}, p^{\lambda - a_1}, 
\dots , p^{\lambda - a_g})$. By definition,
$$S_G(f)(X)= \beta(X) \int _{U(\n Q_p)} f(Xu) du\eqno{(1.1.1)}$$
where $\beta(X)= p^{-ga_1 - (g-1) a_2 - \dots -a_g}$ and $U= 
\left\{\left( \matrix 1 & * \\ 0 & 1 \endmatrix \right)\right\} $, 
entries are $g \times g$-matrices (the multiplier $\beta(X)$ 
differs slightly from the one of [FCh]). 
\medskip
{\bf 1.2.} Here we apply (1.1.1) to $\tau_p = G(\n Z_p)T_pG(\n Z_p)$. Let $f$ be its the characteristic function, and $X=F_i$. 

For $u = \left( \matrix    E_g & A \\ 0 & E_g \endmatrix \right) $ where $A=\left( \matrix u_{11} & u_{12}\\ u_{12}^t & u_{22} \endmatrix \right) $, sizes of diagonal blocks here and below are $g - i$, $i$, we have $F_iu= \left( \matrix B&C\\0&D \endmatrix \right) $, where $B= \diag (1,\dots,1, p,\dots,p)$, $D= \diag(p,\dots,p, 1,\dots,1)$, $C=\left(\matrix u_{11}& u_{12}\\pu_{12}^t & pu_{22}\endmatrix\right)$. 
Hence, $f(F_i u)=1 \iff $ 
entries of $u_{11}, u_{12} \in \n Z_p$, entries of $u_{22} 
\in {1 \over p} \n Z_p$. This implies that 
$$\int _{U(\n Q_p)} f(F_i u) du = p^{{i(i+1) \over 2}}$$
and $S_G(f)(\Phi_i)=1$. For other $X$ it is easy to see that $\int f(Xu)=0$, i.e. 
$$ S_G(\tau_p) = \Phi_0 + \Phi_1 + \dots + \Phi_g\eqno{(1.2.1)}$$

\medskip
{\bf 1.3.} Here we apply (1.1.1) to $\tau_{p,i} = G(\n Z_p)T_{p,i}G(\n Z_p)$, $i\ge1$. Let $f$ be its the characteristic function, and $X=F_j F_k$. We have 

$F_j F_k = \diag (1,\dots,1, p,\dots,p, p^2,\dots,p^2, p^2,\dots,p^2, p,\dots,p, 1,\dots,1)$, $k > j$, sizes of diagonal blocks here and below are $g-k$, $k-j$, $j$, $g-k$, $k-j$, $j$. 

For $u = \left( \matrix    E_g & A \\ 0 & E_g \endmatrix \right) $ where $A=\left( \matrix 
    u_{11} & u_{12} & u_{13}
\\  u_{12}^t & u_{22}& u_{23}
\\  u_{13}^t & u_{23}^t& u_{33} \endmatrix \right) $ we have $F_j F_k u= \left( \matrix B&C\\0&D \endmatrix \right) $, where $B= \diag (1,\dots,1, p,\dots,p, p^2,\dots,p^2$), 

$D= \diag(p^2,\dots,p^2, p,\dots,p, 1,\dots,1)$, $C=\left( \matrix u_{11} & u_{12} & u_{13}
\\ pu_{12}^t & pu_{22}& pu_{23}
\\ p^2u_{13}^t & p^2u_{23}^t& p^2u_{33}
\endmatrix \right) $. Hence, $f(F_j F_k u)=1 \iff $ 
entries of $u_{11}, u_{12}, u_{13} \in \n Z_p$, entries of 
$u_{22}, u_{23} \in {1 \over p} \n Z_p$, entries of $u_{33} 
\in {1 \over p^2} \n Z_p$, $\rank (\widetilde {pu_{22}}) 
= k-j-i$, where tilde means the residue map $\n Z_p \to 
\n F_p$. 
(This is because for a symmetric $g \times g$-matrix $A$ such that $\rank \tilde A=r$ we have $\left( \matrix p & A \\ 0 & p \endmatrix \right) \in 
G(\n Z_p)T_{p,g-r}G(\n Z_p)$). 

So, we denote by $R_g(i) = R_g(i,p)$ the quantity of symmetric 
$g \times g$-matrices with entries in $\n F_p$ of corank exactly 
$i$ (see [AZh], Chapter 3, Lemma 6.19 for the formula for $R_g(i)$) and we have 
$$\int _{U(\n Q_p)} f(F_j F_k u) du = 
R_{k-j}(i)\cdot p^{j(k-j)+j(j+1)}$$
and 
$$S_G(\tau_{p,i})(F_j F_k)= \beta(F_j F_k)
\int _{U(\n Q_p)} f(F_j F_k u) du =
R_{k-j}(i)\cdot 
p^{-\left( \matrix k-j+1 \\ 2 \endmatrix \right)}$$
For other $X$ it is easy to see that $\int f(Xu)=0$, i.e. we have ($i\ge 1$):
$$S_G(\tau_{p,i})=\sum_{j,k\ge0, j+i \le k}^g
R_{k-j}(i)\cdot 
p^{-\left( \matrix k-j+1 \\ 2 \endmatrix \right)}\Phi_j \Phi_k\eqno{(1.3.1)}$$
\medskip
{\bf Remark 1.4.} The above formulas can be used for finding of the Hecke polynomial of $X$. Any element of $\Bbb H(G)$ defines a correspondence on $X$. We denote the algebra of these correspondences by $\n T_p$, it is the quotient ring of $\Bbb H(G)$ by the only relation $\tau_{p,g} = \id$. 

Let us consider the (good) reduction of $X$ at $p$, denoted by $\tilde X$. We denote by $\Corr (\tilde X)$ its algebra of correspondences. Obviously there exists an inclusion $\gamma : \n T_p \to \Corr (\tilde X)$. It is known that it can be included in the commutative diagram: 
$$\matrix S_G: & \n H(G) & \to & \n H(M_s) \\
                     & \beta_1\downarrow & & \beta_2\downarrow \\
				\gamma : & \n T_p & \to & \Corr (\tilde X)
				\endmatrix $$
where $\beta_1$ is the natural projection, $\beta_2$ is an epimorphism with the same kernel $\tau_{p,g} - \id$. 

There is the Frobenius map $f: \tilde X \to \tilde X$, we can 
consider it as a correspondence, i.e. $f \in \Corr (\tilde X)$. 
We have: $f=\beta_2(\Phi_0)$ in $\Corr (\tilde X)$, and $\beta_2(\Phi_g)$ is the Verschiebung correspondence.
The minimal polynomial satisfied by $f$ over 
$\n T_p$ is called the Hecke polynomial.

An explicit algorithm of finding of the Hecke polynomial is a by-product of the calculations of the present paper. There are 2 methods of finding of this polynomial: the first one is to eliminate formally $\Phi_1, \dots , \Phi_g$ from (1.2.1), (1.3.1) and to use the relation $\tau_{p,g}=1$. The second one is to use a description of Langlands parameters of unramified representations --- this gives us formula (2.7.2). See Appendix, Table 4 for the explicit formulas for the cases $g=2,3$. 
\medskip
{\bf Remark 1.5.} There is a slightly another method of the explicit finding of Hecke polynomial. We denote by $\Omega(G)$ the Weyl group of $G$. It enters in the exact sequence $$0\to(\n Z/2\Bbb Z )^g\to\Omega(G)\to S(g)\to 0$$
and there exists a section $i: S(g)\to\Omega(G)$. Let $U_i$, $V_i$ ($i=1,\dots,g$) be independent variables. We have: (see [FCh], Ch. 7 for example) $\Bbb H(T)$ is a subring of $\Bbb Q[U_i^{\pm 1}, V_i^{\pm 1}]$ generated by $(U_iV_i^{-1})^{\pm 1}$, $\prod_{i=1}^gU_i$. $\Omega(G)$ acts on $\Bbb H(T)$ by the obvious manner ($S(g)$ permutes indices in $U_i$, $V_i$, and $(\Bbb Z/2\Bbb Z)^g$ interchanges $U$, $V$). 
Then $\Bbb H(G)$, $\Bbb H(M_s)$ are subrings of $\Bbb H(T)$ stable with respect to $\Omega(G)$, $i(S(g))$ respectively, and Satake maps $S_G$, $S_T$ are identical inclusions. 

For a subset $I$ of $1,\dots,g$ we denote $U_I=\prod_{i\in I}U_i\prod_{i\not\in I}V_i\in \n H(T)$. 
Then we have: $$S_T(\Phi_i)=\sum_{\#(I)=i}U_I\eqno{(1.5.1)}$$ (particularly, $\prod_{i=1}^gV_i$ is the Frobenius element and $\prod_{i=1}^gU_i$ is the Verschiebung). Using (1.2.1), (1.3.1) and (1.5.1) it is easy to find images of $\tau_p$, $\tau_{p,i}$ in $\Bbb Q[U_i^{\pm 1}, V_i^{\pm 1}]$ (for example, $\tau_p=\sum_I U_I= \prod_{i=1}^g (U_i+V_i)$). 

Roots of Hecke polynomial are $(\n Z/2\Bbb Z )^g$-conjugates of $\prod_{i=1}^gV_i$, i.e. elements $U_I$. We denote the $i$-th coefficient of the Hecke polynomial by $\goth h_i\in \Bbb H(G)$. Hence, $\goth h_i=(-1)^i\sigma_i(U_I)$, $i=0,\dots,2^g$, where $\sigma_i$ is the $i$-th symmetric polynomial. $\goth h_i$ can be found explicitly using (1.2.1), (1.3.1).
\medskip
{\bf 2. Description of Weil numbers of $\Cal M_P$.}\nopagebreak
\medskip
{\bf 2.1.} Let $T \subset B \subset G$ be the standard Borel pair, i.e. $T$ is as above and 

$B =\{ \left( \matrix (D^t)^{-1} & * \\ 0 & D \endmatrix 
\right)\in G | D$ is an upper-triangular $g \times g$-matrix. $\}$

Let $\chi : T(\n Q_p) \to \n C^*$ be a nonramified character. 
$\chi$ is defined uniquely by the numbers $$a_0 = \chi( \left( \matrix 1 & 0 \\ 0 & p 
\endmatrix \right) ), a_i = \chi( \left( 
\matrix \rho_i & 0 \\ 0 & \rho_i^{-1} \endmatrix \right) )$$ where 
$\rho_i=\diag (1, \dots , p, \dots , 1)$, $p$ being at the 
$i$-th place, $i= 1, \dots , g$. It is convenient to denote $b_i=p^ia_i$. 

From here and until (2.7) we shall assume that $\chi$ is arbitrary, i.e. $b_i$ are arbitrary numbers. From (2.7) we shall treat only one $\chi$ defined in Introduction, Step 3. 

We can expand $\chi$ on $B(\n Q_p)$ 
using the projection $B \to T$, and let $\pi_{\chi}: G(\n Q_p) \to GL(\Cal V)$ be the parabolically 
induced representation. Recall its definition: $\Cal V$ is a space of functions $f: G(\n Q_p) \to \n C$ which satisfy $$\forall b \in B(\n Q_p) \;\;\; f(bg)=\chi(b)\cdot f(g)$$
and the action is right translation:
$$[\pi_{\chi}(t)(f)](g)=f(gt)$$
There exist a 1-dimensional subspace $\Cal V^{G(\n Z_p)}
\subset \Cal V$ of $G(\n Z_p)$-invariant functions, an action of $\n H(G)$ on 
$\Cal V^{G(\n Z_p)}$ and hence a homomorphism $\alpha_G(\chi): \n H(G) \to \n C$. 

There are 2 methods of description of $\alpha_G(\chi)$: the first one is based on consideration of decomposition of a double 
coset $G(\n Z_p)TG(\n Z_p)$, $T\in G$, as a union of ordinary cosets. Really, if  
$G(\n Z_p)TG(\n Z_p)=\cup_i \gamma_i G(\n Z_p)$ then 
$\alpha_G(\chi)(G(\n Z_p)T G(\n Z_p))= \sum_i \chi(\gamma_i)$. We treat this decomposition in Sections 2.2 - 2.4. 

The second method (which is much more convenient) is treated in 2.5. So, Sections 2.2 -2.4 are entirely of survey nature. 
\medskip
{\bf 2.2.} Here we consider 
for simplicity the case of $G=GL_n$ and a double coset 
$G(\n Z_p)T_{p,i}G(\n Z_p)$ for $T_{p,i}=\diag(1,\dots,1,p,\dots,p)$, $p$ occurs $i$ times. 
This coset decomposition is the following: 
$$G(\n Z_p)T_{p,i}G(\n Z_p)=\bigcup_I \bigcup_{\{c_{jk}\} }\gamma_{I,\{c_{jk}\} } G(\n Z_p)$$
where $I$ runs through the set of all subsets of $\{1,\dots,n\}$ 
containing $i$ elements, $c_{jk}$ belongs to a fixed set of 
representatives of $\n F_p$ in $\n Z$, $c_{jk}=0$ unless $j \not\in I$,
$k \in I$, $j<k$, and 
$$\gamma_{I,\{c_{jk}\} }= \sum_{j \in I} p\cdot e_{jj} + 
\sum_{j \not\in I} e_{jj} +
\sum_{j , k} c_{jk} \cdot e_{kj} $$
($j \not\in I$, $k \in I$, $j<k$), where $e_{lm}$ are elementary matrices. 
\medskip
We can transform the above decomposition as follows: 
$$G(\n Z_p)T_{p,i}^{-1}G(\n Z_p)=\bigcup_* G(\n Z_p) \gamma_*^{-1} = 
\bigcup_*  {\gamma_*^{-1}}^{t} G(\n Z_p);$$
$$G(\n Z_p)pT_{p,i}^{-1}G(\n Z_p)=G(\n Z_p)T_{p,n-i}G(\n Z_p)= 
\bigcup_*  p{\gamma_*^{-1}}^{t} G(\n Z_p)$$
So, we have:
$$\gamma_{I,\{c_{jk}\} }^{-1}= \sum_{j \in I} p^{-1}\cdot e_{jj} + 
\sum_{j \not\in I} e_{jj} +
\sum_{j , k} -p^{-1}c_{jk} \cdot e_{kj} $$
($j \not\in I$, $k \in I$, $j<k$)
and hence
$$p{\gamma_*^{-1}}^{t} = \sum_{j \in I} e_{jj} + 
\sum_{j \not\in I} pe_{jj} +
\sum_{j , k} -c_{jk} \cdot e_{jk} $$
($j \not\in I$, $k \in I$, $j<k$). These elements are in $B$. Further, 
for a fixed $I$ we have $$\chi(p{\gamma_{I,\{c_{jk}\} }^{-1}}^{t})=
\prod_{i\not\in I}a_i$$
and hence 
$$\alpha_G(\chi)(T_{p,n-i})=\sum_{I,\{c_{jk}\} }
\chi(p{\gamma_{I,\{c_{jk}\} }^{-1}}^{t})=
\sum_{I,\#I=i}\prod_{i\not\in I}a_i \cdot p^{\#\{(j,k)\vert 
j \not\in I, k \in I, j<k\} }$$
which gives us $$\alpha_G(\chi)(T_{p,n-i})=p^{-\frac{i(i+1)}{2}}
\sigma_{i}(b_*)$$
\medskip
{\bf 2.3.} Here we consider the case $G=GSp_{2g}(\n Q)$, $T=T_p$. We have the 
following decomposition: $G(\n Z_p)T_pG(\n Z_p) = \cup_i G(\n Z_p)\gamma_i$ where the set 
$\{\gamma_i \}$ is described as follows: 

1. We consider all subsets $I \subset \{1, \dots, g\}$ (there are $2^g$
of them); 

2. If such $I$ is fixed then we consider the set of $\gamma = \left( \matrix A & B \\ 0 & D 
\endmatrix \right) $ such that  
$$D= \sum_{j \in I} p\cdot e_{jj} + 
\sum_{j \not\in I} e_{jj} +
\sum_{j , k} c_{jk} \cdot e_{jk},$$ 
$$A= p{D^t}^{-1}= \sum_{j \in I} e_{jj} + p\cdot
\sum_{j \not\in I} e_{jj} +
\sum_{j , k} -c_{jk} \cdot e_{kj},$$ ($j \not\in I$, $k \in I$, $j<k$)
$$B = \sum_{j,k \in I} b_{jk}e_{jk},$$
$b_{jk}$, $c_{jk}$ belong to a fixed set of 
representatives of $\n F_p$ in $\n Z$, and $b_{jk}=b_{kj}$. 

Now we use the same transformations as above. We have: 
$\chi(p \gamma_i^{-1})= a_0 \prod_{i \in I} a_i$ and it is easy to see that 
$$\alpha_G(\chi)(T_p)=a_0\prod_{i=1}^g(1+b_i)\eqno{(2.3.1)}$$
\medskip
{\bf 2.4.} Here we consider the case $G=GSp_{2g}(\n Q)$, $T=T_{p,i}$. 
Firstly we describe a set $J$ such that 
$$\bigcup_{i=0}^g G(\n Z_p)T_{p,i}G(\n Z_p) = \cup_{j\in J} G(\n Z_p)\gamma_j$$
and then for each $j \in J$ we find the corresponding $i \in 0,\dots, g$. 

We have: $\gamma_j= \left( \matrix A & B \\ 0 & D 
\endmatrix \right) \in GSp_{2g}$ with $\lambda(\gamma_j)=p^2$. $D$ is 
an upper-triangilar matrix whose diagonal entries $D_{ii}$ are $p^{d_i}$, 
$d_i = 0,1,2$, i.e. we have $3^g$ possibilities for the choice of $d_i$. 
To choose a set $d_i$ is the same as to choose a partition $\{1,\dots,g\}
=I_0 \cup I_1 \cup I_2$, $i \in I_k \iff d_i=k$. Non-diagonal entries of $D$
are described as follows: 

(1) If $i \in I_0, j \in I_1, i < j $ then $D_{ij} $ runs through 
a system of representatives in $\n Z$ of $\n Z/p$; 

(2) If $i \in I_0, j \in I_2, i < j $ then $D_{ij} $ runs through 
a system of representatives in $\n Z$ of $\n Z/p^2$; 

(3) If $i \in I_1, j \in I_1, i < j $ then $D_{ij} $ runs through 
a system of representatives in $\n Z$ of $\n Z/p$, and 
the Jordan normal form of this part of $D$ has blocks of size
1 or 2 (i.e. its square is 0); 

(4) If $i \in I_1, j \in I_2, i < j $ then $D_{ij}=pD'_{ij}$, 
where $D'_{ij}$ runs through 
a system of representatives in $\n Z$ of $\n Z/p$; 

Other $D_{ij}$ are 0. We denote submatrices of $D$ described in (1) - (4) 
above by $\goth A$, $\goth B$, $\goth C$, $p\goth D$ respectively.
\medskip
Further, we have $A=p^2D^{-1t}$, and the description of 
$B=\{B_{ij}\}$ is the following. Firstly, 
$B_{ij}=0$ if $i\in I_0$ or $j\in I_0$. Further, we denote 
submatrices of $B$ formed by elements $B_{ij}$ with 
$i\in I_r$, $j\in I_s$ ($r,s=1,2$) by $\goth B_{rs}$. Entries of $\goth B_{11}$, $\goth B_{21}$
(resp. $\goth B_{12}$, $\goth B_{22}$) run through 
a system of representatives in $\n Z$ of $\n Z/p$, (resp. of $\n Z/p^2$). 

Finally, the above matrices satisfy the following relations (which are 
equivalent to a condition $\gamma_j \in GSp_{2g}$): 

\medskip
(1) $\goth B_{11}^t(pI + \goth C) = (pI + \goth C^t) \goth B_{11}$

(2) $(pI + \goth C^t) \goth B_{12} = p\goth B_{11}^t \goth D + 
p^2 \goth B_{21}^t$

(3) $\goth D^t \goth B_{12} + p \goth B_{22} = 
\goth B_{12}^t \goth D + p \goth B_{22}^t$
\medskip
For a given $\gamma_j$ it is possible to find $i$ such that 
$\gamma_j \in G(\n Z_p)T_{p,i}G(\n Z_p)$. It is obvious that $i \le \#I_1$.

For each set ${\goth d}=\{d_i\}$
we denote by $C({\goth d},k)$ the quantity of 
matrices $\gamma_j$ described above such that $\gamma_j \in G(\n Z_p)T_{p,k}G(\n Z_p)$.
In these notations we have the following formula: 
$$\alpha_G(\chi)(T_{p,k})= \sum_{\goth d} C(\goth d,k) 
\prod_{i=1}^g a_i^{d_i}$$ 
Really, it is more convenient to denote $\tilde C({\goth d},k)=C({\goth d},k)p^{-\sum_{i=1}^g id_i}$, so 
$$\alpha_G(\chi)(T_{p,k})= \sum_{\goth d} \tilde C(\goth d,k) 
\prod_{i=1}^g b_i^{d_i}\eqno{(2.4.1)}$$
Formulas for $\tilde C({\goth d},k)$ and $\alpha_G(\chi)(T_{p,k})$ for $g=2,3$ are given in the appendix, tables 5, 6. 
\medskip
{\bf 2.5.} It is well-known that there exists a map $\alpha_T(\chi): \n H(T)\to \n C$ given by the formula 
$$\alpha_T(\chi)(V_i)=a_0^{1/g}, \; \; \alpha_T(\chi)(U_i)=a_0^{1/g}b_i\eqno{(2.5.1)}$$ 
such that 
$$\alpha_G(\chi)=\alpha_T(\chi)\circ S_T\circ S_G\eqno{(2.5.2)}$$
Combining (2.5.1), (2.5.2) with (1.2.1), (1.3.1), (1.5.1), we get 
$$\alpha_G(\chi)(\tau_{p,i})=a_0^2\sum_{j,k\ge0, j+i \le k}^g
R_{k-j}(i)\cdot 
p^{-\left( \matrix k-j+1 \\ 2 \endmatrix \right)}\sum_{\#(J)=j}b_J\sum_{\#(K)=k}b_K\eqno{(2.5.3)}$$
Comparing (2.4.1) and (2.5.3) we get immediately that for $i=1,\dots,g$
$$\tilde C(\goth d, i)=\sum_{j=0}^{[(q_1-1)/2]}R_{q_1-2j}(i)p^{-(q_1-2j+1)(q_1-2j)/2}\left( \matrix q_1\\j\endmatrix \right)\eqno{(2.5.4)}$$
where $q_1=\#I_1$ is the quantity of ones in $\goth d$ and $[x]$ is the integer part of $x$. 
\medskip
{\bf 2.6.} Here we recall a description of the finite-dimensional representation $r:$ $^L G \to GL(W)$ ([BR], 5.1), and its restriction to $\hat T \subset ^L G$ for our case $G=GSp_{2g}(\n Q)$ (this is well-known, see for example [BR], 5.1, Example C). So, firstly we describe the spin representation and its restriction on Cartan subalgebra. The below facts can be found in many sources; we use [J]. 
\medskip
The dual of $GSp_{2g}$ is the spinor group $GSpin_{2g+1}$. Since $G=GSp_{2g}$ is over $\Bbb Q$, we have: $^L G=W_{\Bbb Q}\times GSpin_{2g+1}$, and $r: ^L G \to GL(W)$ (see, for example, [BR], (5.1)) is trivial on $W_{\Bbb Q}$. It is known that $r: GSpin_{2g+1}\to GL(W)$ is the spin representation. There exists a 2-fold covering $\eta: GSpin_{2g+1} \to GO_{2g+1}$. Recall the definition of the corresponding representation of Lie algebras $\goth r: \goth{GO}_{2g+1}\to GL(W)$. 

Let $V$ be a vector space of dimension $2g+1$, $u_1, \dots ,
u_{2g+1}$ its basis and $B$ a quadratic form whose matrix 
in this basis is $\left( \matrix 1&0&0 \\ 0&0&1 \\ 0&1&0 
\endmatrix \right)$, the size of diagonal entries is 
$1,g,g$. We consider the corresponding orthogonal Lie algebras $\goth {GO}(B)$, $\goth {O}(B)$. 
Their Cartan subalgebras of diagonal matrices $\goth {TGO}(B)$, resp. $\goth {TO}(B)$ have bases $\theta_0, \theta_1, \dots , \theta_g$, resp. $\theta_1, \dots , \theta_g$, where $\theta_0$ is the $2g+1\times 2g+1$ unit matrix and $\theta_i = e_{i+1,i+1}-e_{i+g+1,i+g+1}$ for $i>0$,  $e_{ij}$ is an elementary matrix ([J], p. 139, (63)). 

The Clifford algebra $C=C(V,B)$ is the quotient 
of $\sum_{n=0}^{\infty} V^{\otimes n}$ 
(the tensor algebra of $V$) by relations $v_1 \otimes v_2 +
v_2 \otimes v_1 = 2B(v_1,v_2)$. Let $L(C)$ be the 
corresponding Lie algebra, $M_1 \subset C$ the 
natural projection of $ V=V^{\otimes 1} \subset \sum_{n=0}^{\infty} V^{\otimes n}$ to $C$, and $M_2=[M_1, M_1]$. It is known ([J], p. 231,  Th. 7) 
that $M_2$ is a Lie subalgebra of $L(C)$, and it is 
isomorphic to $\goth O(B)$. Further, $M_1$ is isomorphic to $V$ as 
a vector space, and the Lie action of $M_2$ on $M_1$ defined 
by the formula $x(y)=xy-yx$ (here $x \in M_2$, $y \in M_1$, 
multiplication is in $C$), coincides with the action (of a 
matrix on a vector) of $\goth O(B)$ on $V$. 

This formula permits to get an explicit identification of $\goth O(B)$ and $M_2$. Namely, we denote $v_i = u_1 u_{i+1}$, $w_i = u_1 u_{i+1+g}$, 
multiplication is in $C$, $v_i$, $w_i \in M_2$. We have: 
$$\hbox{ for $i>0$ }\ \ \ \theta_i \in \goth O(B) \hbox{ corresponds to }
\frac12 + \frac12 v_iw_i \in M_2\eqno{(2.6.1)}$$ 
(calculations are similar to [J], p. 233, (34) or can be deduced easily from these formulas; it is necessary to take into consideration that $h_i$ of page 139 are $\theta_i$ and $h_i$ of page 233 are $\theta_i - \theta_{i+1}$). 

For $I=(\alpha_1, \dots , \alpha_k)\subset (1,\dots,g)$ we set $x_I=v_1 \cdot \dots \cdot v_g \cdot w_{\alpha_1} \cdot \dots \cdot w_{\alpha_k}$. 
The space of spin representation $W$ is a subspace of 
$C$ spanned on all vectors $x_I$. The action of $\goth O(B)$ is the 
right multiplication by the corresponding elements of $M_2$. This is exactly $\goth r$ restricted on $\goth O(B)$. (2.6.1) shows that 

$$\theta_i(x_{I})=\epsilon x_{I}\eqno{(2.6.2)}$$
where $\epsilon =\frac12$ if $i\in I$ and $\epsilon =-\frac12$ if $i\not\in I$. 

Finally, it is known that $\theta_0$ acts on $W$ by multiplication by $\frac12$.

Let $\hat T\subset GSpin_{2g+1}$ be the dual torus of $T\subset GSp_{2g}$, $\hat \goth T$ its Lie algebra and $\eta_{Lie}: \hat \goth T \to \goth {TGO}(B)$ the restriction of $\eta$. For $t \in T$ we set $t= \diag (x_1, \dots, x_g, \lambda x_1^{-1}, \dots, \lambda x_g^{-1})$, so we can consider $\lambda, x_1, \dots, x_g$ as a basis of $X^*(T)$. We denote the dual basis of $X^*(\hat T)$ by $\lambda', x_1', \dots, x_g'$ and we consider $\lambda', x_1', \dots, x_g'$ as coordinates of an element $t \in \hat T$. Further, we denote by $\nu_0,\nu_1, \dots, \nu_g$ the basis of $\hat \goth T$ dual to $\lambda', x_1', \dots, x_g'$. Formulas for $\eta_{Lie}$ in bases $\nu_*$, $\theta_*$ are the following: 

$$\eta_{Lie}(\nu_0)=2\theta_0, \;\;\;\; \eta_{Lie}(\nu_i)=\theta_i+ \theta_0 \eqno{(2.6.3)}$$
(2.6.2), (2.6.3) imply formulas for the action of $\nu_i$ on $x_I$: $\goth r(\nu_i)(x_I)=x_I$ if $i=0$ or $i \in I$, $\goth r(\nu_i)(x_I)=0$ if $i \not\in I$. In its turn, these formulas imply formulas for $r|_{\hat T}$: $$r(\lambda', x_1', \dots, x_g' )(x_I)=(\lambda'\prod_{i\in I}x_i')x_I\eqno{(2.6.4)}$$
(see also [BR], end of (5.1)). 
\medskip
{\bf 2.7.} Here we apply $r$ to a Langlands element $\theta_{\pi_{\chi}}\in \; ^L G$ of $\pi_{\chi}$ in order to find Weil numbers of $\Cal M_P$. 

We can choose $\theta_{\pi_{\chi}}$ in $\hat T$; namely, it is known that
$\lambda', x_1', \dots, x_g'$-coordinates of $\theta_{\pi_{\chi}}$ are  $$(a_0,b_1,\dots,b_g)\eqno{(2.7.1)}$$
(2.6.4) and (2.7.1) imply that 
\medskip
(2.7.2) $\forall I\subset \{1,\dots,g\}$ the element $x_I$ is an eigenelement of $r(\theta_{\pi_{\chi}})$ with eigenvalue $a_0b_I=\alpha_T(\chi)(U_I)$. 
\medskip
From here we fix $\chi$ such that $\pi_{\chi}=\pi_p$ --- the $p$-part of $\pi$ of Introduction. (2.7.2) and (0.1) give us immediately 
\medskip
(2.7.3) The $2^g$ Weil numbers of $\Cal M_{P}$ have the form $a_0b_I$. 
\medskip
Moreover, the existence of pairing in cohomology of $X$ shows that numbers $b_i$ satisfy a relation $a_0^2\prod_{i=1}^g b_i=p^{g(g+1)/2}$ ($\iff a_0^2\prod_{i=1}^g a_i=1)$. 

Since $\alpha_G(\chi)=\goth m$ ($\goth m$ of the Introduction), (2.3.1), (2.5.3) and (2.7.3) give us expressions of $\goth m(\tau_{p,*})$ in terms of Weil numbers of $\Cal M_P$: 

$$\goth m(\tau_p)=\sum_{I\in 2^g}a_0b_I\eqno{(2.7.4a)}$$ $$\goth m(\tau_{p,i})=a_0^2\sum_{j,k\ge0, j+i \le k}^g R_{k-j}(i)\cdot p^{-\left( \matrix k-j+1 \\ 2 \endmatrix \right)}\sum_{\#(J)=j}b_J\sum_{\#(K)=k}b_K\eqno{(2.7.4b)}$$ where $b_i$ should be interpreted as numbers entering in the formula (2.7.3) for Weil numbers of $\Cal M_P$. 
\medskip
{\bf Remark.} The above formulas give us a simple proof that $\alpha_G(\chi)(\sum \goth h_iT^i)$ is the characteristic polynomial of $r(\theta_{\pi_{\chi}})$ (this is well-known; see, for example, [B] for less explicit proof in a more general situation). Really, roots of the Hecke polynomial are $U_I$ ($I$ runs over $2^g$), and $\alpha_T(\chi)(U_I)=r_I(\theta_{\pi_{\chi}})$. 
\medskip
\medskip
{\bf 3. Some explicit formulas for archimedean cohomological representations of $G$.}\nopagebreak
\medskip
{\bf 3.1. Description of parabolic subgroups of $G$.}\nopagebreak 
\medskip
The set of simple positive roots that corresponds to a Borel pair $(T, B)$ of $G$ is: 
$$\omega_0=x_1^2\lambda^{-1}, \omega_i=x_{i+1}x_i^{-1}, \ \ \ i=1, \dots, g-1,$$
$\lambda, x_i$ of 2.6. We denote this set by $\Delta$. 

Parabolic subgroups that contain $B$ are in one-to-one correspondence to the set of subsets of $\Delta$. We shall tell that such a subgroup is of the first type if the corresponding subset of $\Delta$ does not contain $\omega_0$, and of the second type, if it contains $\omega_0$. The set of subgroups of both types is isomorphic to the set of ordered partitions of $g$, i.e. the set of representations of $g$ as a sum 
$$g=\goth b_1+\goth b_2+\dots+\goth b_k\eqno{(3.1.1)}$$ 
where $\goth b_i\ge 1$, the order is essential. We denote $m_i=\goth b_1+\cdots \goth b_{i-1}$ ($i=1,\dots,k$). The subset of $\Delta$ that corresponds to (3.1.1) is $\Delta-\{\omega_0,\omega_{m_2},\omega_{m_3},\dots,\omega_{m_k}\}$ for the first type, $\Delta-\{\omega_{m_2},\omega_{m_3},\dots,\omega_{m_k}\}$ for the second type. We denote the corresponding parabolic subgroup by $P$ and its Levi decomposition by $P=MN$. Their description is the following: 
\medskip
First type: 
$$M=\left( \matrix A&0 \\ 0&D \endmatrix \right)\eqno{(1M)}$$
where $A$, $D$ are block diagonal matrices with sizes of blocks $\goth b_1,\goth b_2,\dots,\goth b_k$. We denote block entries by $A_i$, $D_i$ respectively $(i=1,\dots,k)$;
$$P=\left( \matrix A&B \\ 0&D \endmatrix \right)\eqno{(1P)}$$
where $A$ (resp. $D$) is a lower (resp. up) block triangular matrix (with the same size of blocks clearly), and 

{\bf (1N)} $N\subset P$ is its subset of matrices whose block entries are unit matrices. 

For the second type we have 
$$M=\left( \matrix A&B \\ C&D \endmatrix \right)\eqno{(2M)}$$
where $A$, $D$ are like in (1M), and $B$, $C$ contain only the upper left corner of size $\goth b_1$ of non-zero elements. These matrices are denoted by $B_1$, $C_1$ respectively; clearly the $2\goth b_1\times2\goth b_1$-matrix $\left( \matrix A_1&B_1 \\ C_1&D_1 \endmatrix \right)\in GSp_{2\goth b_1}$;
$$P=\left( \matrix A&B \\ C&D \endmatrix \right)\eqno{(2P)}$$
where $A$, $D$ are like in (1P), $C$ is like in (2M); 

$$N=\left( \matrix A&B \\ 0&D \endmatrix \right)\eqno{(2N)}$$
where $A$, $D$ are like in (1N), and the upper left corner of size $\goth b_1$ of $B$ is the 0-matrix. 
\medskip
To apply formulas of [BR] we need to describe a Borel pair $(T_c,B_c)$ such that $T_c(\n R)$ 
is compact modulo $Z(\n R)$, where $Z$ is the center of $G$. Namely, 

$T_c$ is the set of matrices $\left( \matrix X&Y \\ -Y&X \endmatrix \right)$ where $X$, $Y$ are diagonal $g\times g$-matrices such that 

$X^2+Y^2=\lambda E_g$.

Let $\alpha=\left( \matrix A&iD \\iA&D \endmatrix \right)$ ($i=\sqrt{-1}$) where $A$, $D$ are any scalar 
$g\times g$-matrices such that $AD=\frac12$; let $A=D=\frac1{\sqrt 2}E_g$. We have $$T_c=\alpha T\alpha^{-1}\eqno{(3.1.2)}$$ We denote $M_c=\alpha M\alpha^{-1}$ and analogically for other objects ($N$, $K$ etc.). 
\medskip
{\bf 3.2. Finding of $\Omega_{\n R}(G)$.}\nopagebreak 
\medskip
Here we recall an explicit description of $\Omega_{\n R}(G)$ which is necessary for finding of $\Pi_P$, see 3.3 below. It is possible to use a fact that it contains a subgroup $\Omega(K_c)$ of index 2, but we give a direct calculation. We denote the normalizer by $\goth N$. There is an isomorphism $\goth N(T_c)/T_c=\Omega(G)$ and a section of sets $\gamma:\Omega(G)\to \goth N(T_c)$. (3.1.2) implies that $\goth N(T_c)=\alpha \goth N(T)\alpha^{-1}$. Let $j=1,\dots, g$, $e_j= (1,\dots,1,-1,1,\dots,1)\in (\n Z/2\n Z)^g \subset \Omega(G)$ ($-1$ is at the $j$-th place). A representative of $e_j$ in $\goth N(T)$ is $\goth w_{j}=\left( \matrix E_g-e_{jj}&ie_{jj} \\ie_{jj}&E_g-e_{jj} \endmatrix \right)$. It commutes with $\alpha$, i.e. we can set $\gamma(e_j)=\goth w_{j}$. Equality  $i\prod_{j=1}^g\goth w_j=\left( \matrix 0&-E_g \\-E_g&0 \endmatrix \right)\in G(\n R)$ shows that a representative of $(-1,\dots,-1)\in (\n Z/2\n Z)^g \subset \Omega(G)$ belongs to $\Omega_{\n R}(G)$. 

Further, for $\goth w\in S(g)\subset \Omega(G)$ we denote by $M_{\goth w}$ the $g\times g$-matrix whose $(j,k)$-th entry is $\delta^{{\goth w}(k)}_j$ (the matrix of permutation). Then we have $\gamma({\goth w})=
\left( \matrix M_{\goth w}&0 \\0&M_{\goth w}\endmatrix \right)$. It belongs to $\goth N(T)$, $\goth N_c(T)$ and commutes with $\alpha$. 

This means that $\Omega_{\n R}(G)$ contains a subgroup $X\subset \Omega (G)$ given by an exact sequence $$0\to\n Z/2\n Z\to X\to S(g)\to 0$$ where $\n Z/2\n Z \subset (\n Z/2\n Z)^g$ is the diagonal embedding. Really, it is possible to show that $X=\Omega_{\n R}(G)$, i.e. elements of $(\n Z/2\n Z)^g$, except the diagonal element, cannot be lifted to $G(\n R)$. 

Finally, for a subset $I$ of $\{1,\dots,g\}$ --- or, the same, an element $I\in (\n Z/2\n Z)^g \subset \Omega(G)$ --- we set $\gamma(I)=\prod_{j\in I}{\goth w}_j$, and we denote this element by ${\goth w}_I$. 
\medskip
\medskip
{\bf 3.3. Finding of $\Pi_P$.}\nopagebreak
\medskip
The members of $\Pi_P$ are parametrized by the double coset space $$\Omega(M_c)\backslash \Omega(G)/\Omega_{\n R}(G)$$ ([BR], 4.2). We have: $\Omega(M_c)=S({\goth b_1})\times \dots\times S({\goth b_k})$ for $P$ of the first type and 
$\Omega(M_c)=\Omega(GSp_{2\goth b_1})\times S({\goth b_2})\times \dots\times S({\goth b_k})$ for $P$ of the second type. The set of representatives of $\Omega(G)/\Omega_{\n R}(G)$ can be chosen as a half of $(\n Z/2\n Z)^g$ (we choose one element in each pair of elements $(a, (-1,\dots,-1)a)$, $a\in (\n Z/2\n Z)^g$). The above groups $\Omega(M_c)$ act on this set of representatives from the left, 
hence the invariant of their action is the quantity of $1, -1$ in the  segments of length $\goth b_1,\goth b_2,\dots,\goth b_k$ (first type); $\goth b_2,\dots,\goth b_k$ (second type) in the whole segment of length $g$. This means that the set $\Pi_P$ coincides with 

First type: the set of sequences of numbers $c_1,\dots,c_k$ where $0\le c_j\le \goth b_j$ factorized by the equivalence relation $c_1,\dots,c_k \sim \goth b_1-c_1,\dots,\goth b_k-c_k$; representatives $w$ of the corresponding double cosets are 
$$w=(\underbrace{1,\dots,1}_{c_1\hbox{ times}},\underbrace{-1,\dots,-1}_{\goth b_1-c_1\hbox{ times}},\dots, \underbrace{1,\dots,1}_{c_k\hbox{ times}},\underbrace{-1,\dots,-1}_{\goth b_k-c_k\hbox{ times}})\in (\n Z/2\n Z)^g \subset \Omega(G)\eqno{(3.3.1)}$$
Second type: the same, but the sequences are $c_2,\dots,c_k$, and 
$$w=(\underbrace{1,\dots,1}_{\goth b_1\hbox{ times}},\underbrace{1,\dots,1}_{c_2\hbox{ times}},\underbrace{-1,\dots,-1}_{\goth b_2-c_2\hbox{ times}},\dots, \underbrace{1,\dots,1}_{c_k\hbox{ times}},\underbrace{-1,\dots,-1}_{\goth b_k-c_k\hbox{ times}})\in (\n Z/2\n Z)^g \subset \Omega(G)\eqno{(3.3.2)}$$
Notation: such a sequence $c_1,\dots,c_k$ or $c_2,\dots,c_k$ is denoted by $\goth c$ and the set of all there sequences by $\goth C$. We denote the set of $w\in \Omega(G)$ of the form (3.3.1), (3.3.2) by $\goth W$, i.e. there is a 1 --- 1 correspondence between $\goth C$ and $\goth W$: $w=w(\goth c)$, $\goth c = \goth c(w)$. The representation $\pi\in\Pi_P$ that corresponds to $\goth c$ is denoted by $\pi_{\goth c}$ or (like in [BR]) by $\pi_w$.
\medskip
\medskip
{\bf 3.4. Finding of $p_w$, $q_w$.}\footnotemark \footnotetext{This section is not logically necessary for the proof of the theorem.} \nopagebreak
\medskip
Numbers $p_w$, $q_w$ are defined in [BR], 4.3; here we use notations of this paper. Firstly we recall the definition of $\Cal P_c^{\pm}$ and find them explicitly. Let $h:\Res_{\n C/\n R}G_m\to G$ be a Deligne map for the Siegel variety. We use the following $h$: for $z=x+iy$  $h(z,\bar z)=\left( \matrix x&y \\-y&x \endmatrix \right)$. Let $i_1: G_m \to \Res_{\n C/\n R}G_m$ be the map $z \to (z,1)$ and $\mu=h\circ i_1$. $\Cal P_c^{\pm}$ are the subspaces of $\goth{GSp}_{2g}$ on which $\ad \mu(t)$ acts by $t^{-1}$ and $t$ respectively (see for example [D] or [BR], 4.3).  
An element of $\goth{GSp}_{2g}$ is a matrix $\left( \matrix A&B \\ C&-A^t+(\lambda-1)E_g \endmatrix \right)$ where 
$B$, $C$ are symmetric. A calculation gives us: 
$$\Cal P_c^+=\left( \matrix C&iC \\ iC&-C \endmatrix \right),
\Cal P_c^-=\left( \matrix C&-iC \\ -iC&-C \endmatrix \right)$$ where $C$ is a symmetric $g\times g$-matrix. 

For $w \in \goth W\subset \Omega(G)$ we have $\gamma (w)={\goth w}_I$  for $I=$ the set of $-1'$s in (3.3.1), (3.3.2); we denote it simply by ${\goth w}$. Further, we denote by $\Cal N$, $\Cal N_{\goth w}$, $\Cal N_c$, $\Cal N_{c\goth w}$ the Lie algebras of $N$, ${\goth w}^{-1}N{\goth w}$, $N_c$, ${\goth w}^{-1}N{c\goth w}$ respectively. Numbers $p_w=\dim(\Cal N_{c\goth w}\cap\Cal P_c^+)$, $q_w=\dim(\Cal N_{c\goth w}\cap\Cal P_c^-)$ are defined in [BR], 4.3. It is more convenient to conjugate with $\alpha$: we set $\Cal P^{\pm}=\alpha^{-1} \Cal P_c^{\pm} \alpha$. A calculation gives: 
$\Cal P^{+}=\left( \matrix 0&C \\0&0\endmatrix \right)$,  
$\Cal P^{-}=\left( \matrix 0&0 \\C&0\endmatrix \right)$ where $C$ is a symmetric $g\times g$-matrix. So, 
$p_w=\dim(\Cal N\cap {\goth w}\Cal P^+{\goth w}^{-1})$, $q_w=\dim(\Cal N\cap {\goth w}\Cal P^-{\goth w}^{-1})$. 

Further, $\Cal N$ has the same description like in (1N), (2N), but the  diagonal blocks are 0-matrices. 

Let $e_{i,j}$ be the elementary $(i,j)$-matrix. Matrices ${\goth w}e_{i,g+j}{{\goth w}}^{-1}$ are given by the following table (here and below we indicate in the third column of the table whether ${\goth w}e_{i,g+j}{{\goth w}}^{-1}\in\Cal N$ or not). 
\medskip
First type:\nopagebreak
\medskip
\settabs 4 \columns
\+ Subtype &${\goth w}e_{i,g+j}{{\goth w}}^{-1}$ \cr \nopagebreak
\medskip
\+ 1. $i\not\in I$, $j\not\in I$&$e_{i,g+j}$& always $\in\Cal N$ \cr \nopagebreak
\+ 2. $i\in I$, $j\not\in I$&$e_{g+i,g+j}$& $\in\Cal N\iff j>i$ and (*)  \cr \nopagebreak
\+ 3. $i\not\in I$, $j\in I$&$e_{i,j}$& $\in\Cal N\iff i>j$ and (*)\cr \nopagebreak
\+ 4. $i\in I$, $j\in I$&$e_{g+i,j}$& never $\in\Cal N$ \cr

where (*) means: $i$, $j$ do not belong to the same segment of partition $g=\goth b_1+\cdots+\goth b_k$. 

Since $C$ is a symmetric matrix, we can take always $j\ge i$, and hence 
the quantity of pairs $(i,j)$ such that ${\goth w}e_{i,g+j}{{\goth w}}^{-1}\in\Cal N$ is: 

Subtype 1. $\frac{(g-\sum_{l=1}^k  c_l)(g+1-\sum_{l=1}^k  c_l)}2$;

Subtype 2. $c_1(\goth b_2-c_2+\goth b_3-c_3+\cdots+\goth b_k-c_k)+c_2(\goth b_3-c_3+\cdots+\goth b_k-c_k)+
\cdots+c_{k-1}(\goth b_k-c_k)$, hence

$$p_w=\frac{(g-\sum_{l=1}^k c_l)(g+1-\sum_{l=1}^k  c_l)}2+\sum_{1\le i<j\le k}c_i\goth b_j-\sigma_2(c_*)$$

Analogously, in order to find $q_w$, we have: 
\medskip
\settabs 4 \columns
\+ Subtype &${\goth w}e_{g+i,j}{{\goth w}}^{-1}$\cr \nopagebreak
\medskip
\+ 1. $i\not\in I$, $j\not\in I$&$e_{g+i,j}$& never $\in\Cal N$ \cr \nopagebreak
\+ 2. $i\in I$, $j\not\in I$&$e_{i,j}$& $\in\Cal N\iff i>j$ and (*)\cr \nopagebreak
\+ 3. $i\not\in I$, $j\in I$&$e_{g+i,g+j}$& $\in\Cal N\iff j>i$ and (*) \cr \nopagebreak
\+ 4. $i\in I$, $j\in I$&$e_{i,g+j}$& always $\in\Cal N$ \cr
with the same notations and assumptions, hence 
$$q_w=\frac{(\sum_{l=1}^k  c_l)(1+\sum_{l=1}^k  c_l)}2+\sum_{1\le j<i\le k}c_i\goth b_j-\sigma_2(c_*)$$

Type 2 is analogous to the type 1. We set $c_1=0$, the above tables are the same with the 
following exception: for subtype 1 (i.e. $i\not\in I$, $j\not\in I$) we have: ${\goth w}e_{i,g+j}{{\goth w}}^{-1}=e_{i,g+j}\in\Cal N$ always except $i,j\in [1,\goth b_1]$. This changes the value of $p_w$:

$$p_w\hbox{ of the second type }=p_w\hbox{ of the first type }-\frac{\goth b_1(\goth b_1+1)}2$$
$$=\frac{(g-\sum_{l=2}^k  c_l)(g+1-\sum_{l=2}^k c_l)}2+\sum_{2\le i<j\le k}c_i\goth b_j-\sigma_2(c_*)-\frac{\goth b_1(\goth b_1+1)}2$$
and for $q_w$ we have the same formula like in the first type: 
$$q_w=\frac{(\sum_{l=2}^k  c_l)(1+\sum_{l=2}^k  c_l)}2+\sum_{1\le j<i\le k}c_i\goth b_j-\sigma_2(c_*)$$
\medskip
{\bf Remarks.} 1. Change of $(c_1, \dots, c_k)$ to $(\goth b_1-c_1, \dots, \goth b_k-c_k)$ leads to interchange of $p_w$, $q_w$. 
\medskip
2. We have: $p_w+q_w=\frac{g(g+1)}2-\sum_{l=1}^k  c_l(\goth b_l-c_l)$ (type 1), 

$p_w+q_w=\frac{g(g+1)}2-\frac{\goth b_1(\goth b_1+1)}2-\sum_{l=2}^k  c_l(\goth b_l-c_l)$ (type 2).
\medskip
{\bf 3.5. Finding of length of representatives of $\Omega(M_c)/\Omega(M_c\cap {\goth w}K_c{\goth w}^{-1})$. }
\medskip
We continue to work with the same $w \in \goth W$, $\goth w\in \goth N(T_c)$ from 3.4. To prove proposition 4.3 below, we must find representatives of the minimal length of $\Omega(M_c)/\Omega(M_c\cap {\goth w}K_c{\goth w}^{-1})$, and find their length (see [BR], 4.3 or [A], proof of (9.1)). Firstly we find $K_c$ --- the centralizer of $\mu$ in $G(\n R)$. It is clear that $K_c$ is the centralizer of $h(\Res_{\n C/\n R}G_m)$ as well. Replacing $h$ by $\alpha^{-1} h \alpha$ we see that $\im \alpha^{-1} h \alpha =\left\{\left( \matrix Z&0 \\0&\lambda Z^{-1}\endmatrix \right)\right\}$ where $Z$ is a scalar matrix. 

We define $K$ to be the centralizer of $\im \alpha^{-1} h \alpha$ in $G$; we have: 

1) $K=\left\{\left( \matrix A&0 \\0&\lambda A^{t-1}\endmatrix \right)\right\}$ where $A\in GL_g$; 

2) $K_c=\alpha K \alpha^{-1}$;

3) $\Omega(T_c,K_c)=\Omega(T,K)=S(g)$. 

Now we see that conjugating with $\alpha$ we get $\Omega(M_c)/\Omega(M_c\cap {\goth w}K_c{\goth w}^{-1})=\Omega(M)/\Omega(M\cap {\goth w}K{\goth w}^{-1})$. 
Like in (3.4), we have a table of ${\goth w}$-conjugates of elementary matrix $e_{i,j}$ ($1\le i,j\le g$): 
\settabs 10 \columns
\+ Subtype &&${\goth w}e_{i,j}{{\goth w}}^{-1}$\cr \nopagebreak
\+ 1. $i\not\in I$, $j\not\in I$&&$e_{i,j}$&& $\in M \iff (*)$ \cr \nopagebreak
\+ 2. $i\in I$, $j\not\in I$ &&$e_{g+i,j}$&& $\not\in M$ (Type 1); $\in M \iff i,j\in (1,\dots,\goth b_1)$ (Type 2) \cr \nopagebreak
\+ 3. $i\not\in I$, $j\in I$ &&$e_{i,g+j}$&& $\not\in M$ (Type 1); $\in M \iff i,j\in (1,\dots,\goth b_1)$ (Type 2) \cr \nopagebreak
\+ 4. $i\in I$, $j\in I$ &&$e_{g+i,g+j}$&& $\in M \iff (*)$  \cr
where (*) here means: $ i,j$ belong to the same segment of the partition $g=\goth b_1\cdots+\goth b_k$.

This means that $$M\cap {\goth w}K{\goth w}^{-1}=\prod_{l=1}^gGL(c_l)\times GL(\goth b_l-c_l)\hbox{ (Type 1); }$$ $$=GL(\goth b_1)\times \prod_{l=2}^gGL(c_l)\times GL(\goth b_l-c_l)\hbox{ (Type 2)}$$ is the set of block diagonal simplectic matrices with block sizes $$c_1,\goth b_1-c_1,\dots,c_k,\goth b_k-c_k,c_1,\goth b_1-c_1,\dots,c_k,\goth b_k-c_k \hbox{ (Type 1)}$$; $$\goth b_1,c_2,\goth b_2-c_2,\dots,c_k,\goth b_k-c_k,\goth b_1,c_2,\goth b_2-c_2,\dots,c_k,\goth b_k-c_k \hbox{ (Type 2)}$$, and 
$$\Omega(M_c)/\Omega(M_c\cap {\goth w}K_c{\goth w}^{-1})=\prod_{l=1}^gS(\goth b_l)/(S(c_l)\times S(\goth b_l-c_l))\hbox{ (Type 1); }$$
$$=(\n Z/2\n Z)^{\goth b_1}\times\prod_{l=2}^gS(\goth b_l)/(S(c_l)\times S(\goth b_l-c_l))\hbox{ (Type 2) }$$
The set $S(b)/(S(c)\times S(b-c))$ is isomorphic to $\goth S(c,b)$ --- the set of all subsets of order $c$ of the set $(1,\dots,b)$ (see Introduction, Step 2c). Let $D \in \goth S(c,b)$, $D=(d_1,\dots,d_c)$, where $1\le d_1<\cdots<d_c\le b$. The equivalence class that corresponds to this $D$ is the set of permutations of $(1,\dots,b)$ that send $(1,\dots,c)$ to $(d_1,\dots,d_c)$. Since the length of a permutation (considered as an element of $S(b)=\Omega(GL(b+1))$ ) is the quantity of inversions of elements, it is easy to see that the permutation with the minimal length in the equivalence class corresponding to $D$ is the permutation that sends $j$ to $d_j$ for $j=1,\dots,c$, and analogously (in increasing order) for $j=c+1,\dots,b$. We denote this permutation by $m_D \in S(b)$; 
we have $l(m_D)=\sum_{j=1}^c d_j-\frac {c(c+1)}2$. 

Further, let $\goth a =(\goth a_1,\dots, \goth a_{\goth b_1})\in (\n Z/2\n Z)^{\goth b_1}\subset \Omega(G)$, where $\goth a_i=0,1$. It is known that $l(\goth a)=\sum_{i=1}^{\goth b_1} i\goth a_i$. 

Finally, the set of representatives of $\Omega(M_c)/\Omega(M_c\cap {\goth w}K_c{\goth w}^{-1})$ of minimal length is $\goth S(\goth c, P)$ of (0.3). Really, let $\rho \in \goth S(\goth c, P)$, $\rho=(D_1,\dots,D_k)$ for $P$ of the first type, $\rho=(\goth a,D_2,\dots,D_k)$ for $P$ of the second type, where $D_i$ is a subset of order $c_i$ of the $i$-th segment of the partition $g=\goth b_1\cdots+\goth b_k$ of $(1,\dots,g)$ and $\goth a$ is as above. We have $m_{D_i}\in S(\goth b_i)$. For $P$ of type 1 the representative of minimal length is $m_{\rho} = m_{D_1}\times \dots \times m_{D_k} \in S(\goth b_1) \times \dots \times S(\goth b_k)\subset S(g) \subset \Omega(G)$, and we have 
$$l(m_{\rho})=\sum_{i=1}^kl(m_{D_i})\eqno{(3.5.1)}$$ 
For $P$ of type 2 we let $m'_{\rho} = m_{D_2}\times \dots \times m_{D_k} \in S(\goth b_2) \times \dots \times S(\goth b_k)\subset S(g) \subset \Omega(G)$ and $m_{\rho}=(\goth a $ multiplied semidirectly by  $m'_{\rho})\in\Omega(G)$. We have 
$$l(m_{\rho})=l(\goth a)+\sum_{i=2}^kl(m_{D_i})\eqno{(3.5.2)}$$ 

{\bf Remark.} It is convenient to treat numbers $f_i=d_i-i$ instead of $d_i$, so $f_1\le f_2\le\cdots\le f_c\le b-c$. For the case $P= \left( \matrix * & * \\ 0 & * \endmatrix \right)$ (i.e. $P$ of type 1, $k=1$), $\goth c=\{c\}$, $w=w(\goth c)$ we have $h^{p_w+r,q_w+r}(\goth g,K_c;\pi_w)=$ the quantity of Young diagrams of weight $r$ in the rectangle with sides $c,g-c$. Analogous formulas exist for other $P$. 
\medskip
\medskip
{\bf 4. Relations between Hecke eigenvalues.}\nopagebreak
\medskip
Formulas (3.5.1), (3.5.2) can be used in order to find dimensions of $$H^{i,j}(\goth g,K_c;\pi_w)$$ ([BR], 4.3). They give us also relations between Weil numbers of $\Cal M_P$. The preliminary form of these relations is the following:
\medskip
{\bf Proposition 4.1.} For any $\goth c \in \goth C$ there is a number $\goth x_{\goth c}$ such that the set of all Weil numbers of $\Cal M_P$ is the following: 
$$p^{l(m_{\rho})}\goth x_{\goth c}$$ 
where $\goth c$ runs over $\goth C$, for a fixed $\goth c$  ${\rho}$ runs over $\goth S(\goth c,P)$. 
\medskip
We shall not give a proof of (4.1), because we need a more general proposition 4.3, see below. 
\medskip
Comparing (4.1) with (2.7.2), we get the following problem: 
\medskip
(4.2). Find relations between $b_i=p^ia_i$ such that both (2.7.2), (4.1) are satisfied.
\medskip
The solution of (4.2) --- and even a more exact result --- is given by the following proposition. Recall that $m_j=\goth b_1+\cdots \goth b_{j-1}$ ($j=1,\dots,k$). 
\medskip
{\bf Proposition 4.3.} First type: $b_{m_j+1}$ are free variables, $b_{m_j+i}=p^{i-1}b_{m_j+1}$ ($i=1,\dots,\goth b_j$), and $a_0$ is defined by the equality $a_0^2\prod _{i=1}^g b_i=p^{g(g+1)/2}$. 

Second type: $b_j=p^j$ for $j=1,\dots, \goth b_1$, $b_{m_j+1}$ ($2\le j\le k$) are free variables, $b_{m_j+i}$ and $a_0$ are like in the first type. 
\medskip
{\bf Proof.} It follows immediately from [A], proof of Proposition 9.1. Let us recall some definitions of loc. cit., page 62 (here $\pi\in\Pi_P$, $\goth{g=GSp}_{2g}$): $$V_{\pi}=\bigoplus_i H^i(\goth g,K_c;\pi), \ \ \  V_P=V_{\Psi_P}=\bigoplus _{\pi\in\Pi_P} V_{\pi}$$ ($V_P$ is denoted in [A] by $V_{\psi}$ and is defined on the page 59, two lines below (9.2)).  Spaces $V_{\pi}$, $W$ are $\goth{sl}_2(\n C)$-modules (see loc. cit. for the definition of the action of $\goth{sl}_2(\n C)$), and all $V_{\pi}$ and hence $V_P$ have the Hodge decomposition. 

There exist bases $B(W)$, $B(V_{\pi})$, $B(V_P)$ of $W$, $V_{\pi}$, $V_P$ respectively and an isomorphism $\goth d: B(V_P) \to B(W)$ (see [A], line below (9.6)) which gives an isomorphism of $\goth{sl}_2(\n C)$-modules $V_P \to W$. 

Arthur uses a slightly another description of $B(W)$ than the one used in (2.6). Namely, the set of elements of $B(W)$ is isomorphic to the set of cosets $\Omega (G)/\Omega (K_c)$, where $\Omega (K_c) \overset{i}\to{\hookrightarrow} \Omega (G)$ is equal to $S(g)\overset{i}\to{\hookrightarrow}\Omega (G)$ of Remark 1.5. It is clear that $\Omega (G)/\Omega (K_c) = (\n Z/2\n Z)^g$. Let $I\subset \{1,\dots,g\}$; we can treat $I$ as an element of $(\n Z/2\n Z)^g$ as usually. The element of $B(W)$ that correspond to $I$ according loc.cit. is exactly $x_I$ of (2.6).

Now let $\pi=\pi_w$, $w=w(\goth c)$. The set $B(V_{\pi})$ is isomorphic to $\Omega (M_c)/\Omega (M_c\cap {\goth w}K_c {\goth w}^{-1})=\goth S(\goth c, P)$. For any finite group $A$ and its subgroups $B$, $C$ we have 
$$A/B=\bigcup_{a\in C\backslash A/B} C/(C\cap aBa^{-1})\eqno{(4.3.1)}$$
here and below all unions are disjoint. 

Now we apply (4.3.1) to the case $A=\Omega (G)$, $B=K_c$, $C=\Omega (M_c)$ in order to get an inclusion:
$$B(V_{\pi})=\Omega (M_c)/\Omega (M_c\cap {\goth w}K_c {\goth w}^{-1})=\goth S(\goth c, P) \overset{\goth d_w}\to{\hookrightarrow}\Omega (G)/\Omega (K_c) = (\n Z/2\n Z)^g=B(W)$$ It follows from loc. cit. that for ${\rho} \in \goth S(\goth c, P)$ as in the end of (3.5) we have: $\goth d_w({\rho})=x_I$ where for Type 1: 
$$I=D_1\cup\dots\cup D_k\eqno{(4.3.2)}$$
for Type 2:
$$I=I_{\goth a}\cup D_2\cup\dots\cup D_k\eqno{(4.3.3)}$$
where $I_{\goth a}\subset \{1,\dots,\goth b_1\}$ is the set of ones (additive writing of $(\n Z/2\n Z)^g$) in $\goth a$, the union is in $\{1,\dots,g\}$. The Hodge type of $\goth d_w({\rho})$ is $$p_w+l(m_{\rho}), q_w+l(m_{\rho})\eqno{(4.3.4)}$$

Finally, we have $$B(V_P)=\bigcup_{\pi\in\Pi_{P}}B(V_{\pi})$$ and $\goth d:B(V_P) \to B(W)$ is the union of $\goth d_w$ in the obvious sense. 

Let $X$, $Y$, $H$ be the standard basis of $\goth{sl}_2(\n C)$. We have the following properties of the action of $\ad X$ on $V_{\pi}$ and $W$:
\medskip
{\bf (4.3.5)} If $v\in V_{\pi}$ is of the pure Hodge type $(p,q)$ then $\ad X(v)$ is of the pure Hodge type $(p+1,q+1)$. 
\medskip
{\bf (4.3.6)} If $w\in W$ is a $r(\theta_{\pi_{\chi}})$-eigenelement of eigenvalue $\lambda$, then $\ad X(w)$ is a $r(\theta_{\pi_{\chi}})$-eigenelement of eigenvalue $p\lambda$. 
\medskip
Type 1. We use notations $\goth e_j=(0,\dots,0,1,0,\dots,0)\in \goth C$ (1 is at the $j$-th place), $I(n)$ is a subset of $\{1,\dots,g\}$ consisting of the single element $n$, and we denote $x_{I(n)}$ simply by $x_n$. We fix some $j$ and we set $\goth c=\goth e_j$. (4.3.2) shows that $\goth d(V_{\pi_{\goth c}})$ is generated by $x_{m_{j}+1},\dots,x_{m_{j}+\goth b_j}$. According (4.3.4), $\forall i=1,\dots,\goth b_j$ the Hodge type of $x_{m_{j}+i}$ is $p_{\goth c}+i$, $q_{\goth c}+i$. (4.3.5) implies that 
$$\ad X(x_{m_{j}+i})=c_{j,i} x_{m_{j}+i+1}\eqno{(4.3.7)}$$ where $c_{j,i}$ is some non-0 coefficient. 
Now, (2.7.2), (4.3.6) and (4.3.7) imply immediately that $b_{m_{j}+i+1}=pb_{m_{j}+i}$ which is 4.3 for Type 1. 
\medskip
Type 2. The idea of the proof is the same. Firstly we consider $\goth c=(0,\dots,0)$. $B(V_{\pi_{\goth c}})$ is the set of subsets of $\{1,\dots,\goth b_1\}$. (4.3.3) shows that $\goth d(V_{\pi_{\goth c}})$ is generated by $x_I$, where $I \subset \{1,\dots,\goth b_1\}$. (4.3.5) implies that for $\forall i=1,\dots,\goth b_1$ $$(\ad X)^i(x_{\emptyset})=\sum_{I\subset \{1,\dots,\goth b_1\}\suchthat l(I)=i} c_I x_I\eqno{(4.3.8)}$$
where coefficients $c_I$ can be easily found using methods of [VZ]. For us it is sufficient to use a fact that $c_{I(i)}\ne 0$. (2.7.2), (4.3.6) and (4.3.8) imply by induction by $i$ that $b_i=p^i$. 
\medskip
Finally, we consider $\goth c=\goth e_j$ like in Type 1, but with the first zero omitted. $B(V_{\pi_{\goth c}})=(\n Z/2\n Z)^{\goth b_1}\times \{1,\dots,\goth b_j\}$. (4.3.3) shows that $\goth d(V_{\pi_{\goth c}})$ is generated by $x_{J\cup I(m_{j}+i)}$, where $J \subset \{1,\dots,\goth b_1\}$ and $i\in \{1,\dots,\goth b_j\}$. The Hodge type of $x_{J\cup I(m_{j}+i)}$ is $p_{\goth c}+l(J)+i, q_{\goth c}+l(J)+i$. (4.3.5) implies that for $\forall i=1,\dots,\goth b_j$ $$(\ad X)^{i-1}(x_{\emptyset\cup I(m_{j}+1)})=\sum c_{J,j,n} x_{J\cup I(m_{j}+n)}\eqno{(4.3.9)}$$ the sum is over the pairs $(J,n)$, $J \subset \{1,\dots,\goth b_1\}$, $n\in \{1,\dots,\goth b_j\}$ such that $l(J)+n=i$. 

Again it is sufficient to use a fact that $c_{\emptyset,j,i}\ne 0$. As earlier (2.7.2), (4.3.6) and (4.3.9) imply by induction by $i$ that $b_{m_{j}+i}=p^{i-1}b_{m_{j}+1}$. $\square$
\medskip
{\bf Remark 1.} There are $g-k$ (first type); $g-k+1$ (second type) relations between eigenvalues of $\tau_p$, $\tau_{p,i}$ on $\Cal M_P$. 
\medskip
{\bf Remark 2.} (4.1) is obviously a corollary of (4.3); numbers $\goth x_{\goth c}$ are products of some $b_i$ and powers of $p$. 
\medskip
{\bf Remark 3.} Formulas of (4.3) are not direct corollaries of (4.1), (2.7.2): it is easy to construct an example of numbers $a_i$ having another form as in (4.3) but such that both (4.1), (2.7.2) are satisfied. 
\medskip
 We denote $\goth m(\tau_{p})$, $\goth m(\tau_{p,i})$ by $\goth m_{p}$, $\goth m_{p,i}$ respectively. 
\medskip
{\bf Theorem 4.4.} Relations between $\goth m_{p}^2$, $\goth m_{p,i}$ are linear. Particularly, for the parabolic subgroup $P$ of the second type such that $k=g$, all $\goth b_i$ are 1 (see Appendix, 8b) the only relation between $\goth m_{p}^2$, $\goth m_{p,i}$ is the following:
$$\frac{\goth m_p^2}{(p+1)^2}+\sum_{i=1}^{g}Y_i\goth m_{p,i}=0\eqno{(4.5)}$$ 
where $\goth m_{p,g}=1$ and $Y_i$ are polynomials in $p$ (particularly, they do not depend on $g$) defined as follows: $Y_1=-1$ and $Y_n$ is defined by the recurrence relation
$$[\sum_{i=1}^{n-1}Y_iR_{n-1}(i)](1+p^2)p^{-\frac{(n-1)n}2}+
[\sum_{i=1}^{n}Y_iR_{n}(i)]p^{-\frac{n(n+1)}2+1}$$ $$+
[\sum_{i=1}^{n-2}Y_iR_{n-2}(i)]p^{-\frac{(n-2)(n-1)}2+1}+2=0\eqno{(4.6)}$$ 
\medskip
Proof. Follows immediately from (2.7.4) and (4.3). $\square$ 
\medskip
$Y_2$, $Y_3$ are given in the Appendix, Table 7. 
\medskip
\medskip
{\bf Appendix.}\nopagebreak
\medskip
{\bf 1. Some relations satisfied by $\tau_{p,i}$.}\nopagebreak
\medskip
We set $W_n(p)=\prod_{i=1}^n(p^i+1)$. Let $\deg: \n T_p \to \n Z$ be a map of the degree of a double coset (= the quantity of ordinary cosets in it). We have equalities: $$(\tau_p)^2 = \sum _{i=0}^g   \tau_{p,i} W_i(p); \; \; \; \; \; \;  \deg \tau_p = W_g(p)$$ $$\deg \tau_{p,k} = \sum _{\goth d\in 3^g} C({\goth d},k)$$ $$=p^{(g-k)(g-k+1)/2} \frac{W_g(p)}{W_k(p)} \frac{(p^{g}-1)(p^{g-1}-1)\cdot\dots\cdot(p^{g-k+1}-1)} {(p^{k}-1)(p^{k-1}-1)\cdot\dots\cdot(p^{1}-1)}.$$ Particularly, $\deg \tau_{p,0} =p^{(g)(g+1)/2}W_g(p)$, 
\medskip
for $g=2$: $\deg \tau_{p,0} =p^6+p^5+p^4+p^3$, 

$\deg \tau_{p,1} =p^4+p^3+p^2+p$, 
\medskip
for $g=3$: $\deg \tau_{p,0} =p^{12}+p^{11}+p^{10}+2p^9+p^8+p^7+p^6$, 

$\deg \tau_{p,1}=p^{10}+p^{9}+2p^{8}+2p^7+2p^6+2p^5+p^4+p^3 =p^3(p^2+p+1)(p^2+1)(p^3+1)$, 

$\deg \tau_{p,2} =p^6+p^5+p^4+p^3+p^2+p$
\medskip
{\bf Table 2. Numbers $R_g(k)$.}\nopagebreak
\medskip
Source: [AZh], Chapter 3, Lemma 6.19. We have: $R_g(g)=1$, $R_g(g-1)=p^g-1$.\nopagebreak 
$$\hbox{Numbers }R_g(k): \ \ \ \ \matrix & g & 2& 3\\
          k&&&\\
		  0&&p^3-p^2&p^6-p^5-p^3+p^2\\
		  1&&p^2-1&p^5-p^2\\
		  2&&1&p^3-1 \\
		  3&&&1 \endmatrix$$
\medskip
{\bf Table 3. Explicit formulas for Satake map $S_G$.}\nopagebreak
\medskip
$g=2$:\nopagebreak
$$\tilde \tau_p = \Phi_0 + \Phi_1 + \Phi_2$$
$$\tilde \tau_{p,1} = {1\over p}(\Phi_0 \Phi_1+ 
\Phi_1 \Phi_2)+ {p^2 -1\over p^3}\Phi_0 \Phi_2$$
$$\tilde \tau_{p,2} = {1\over p^3}\Phi_0 \Phi_2$$
\medskip
$g=3$:\nopagebreak
$$\tilde \tau_p = \Phi_0 + \Phi_1 + \Phi_2 + \Phi_3 $$
$$\tilde \tau_{p,1} = {1\over p}(\Phi_0 \Phi_1+ 
\Phi_1 \Phi_2+ \Phi_2 \Phi_3)+ {p^2 -1\over p^3}(\Phi_0 \Phi_2 + 
\Phi_1 \Phi_3) + {p^3 -1\over p^4}\Phi_0 \Phi_3 $$
$$\tilde \tau_{p,2} ={1\over p^3}(\Phi_0 \Phi_2 + 
\Phi_1 \Phi_3) + {p^3 -1\over p^6}\Phi_0 \Phi_3 $$
$$\tilde \tau_{p,3} = {1\over p^6}\Phi_0 \Phi_3 $$
\medskip
{\bf Table 4. Coefficients $\goth h_i$ of the Hecke polynomial.}\nopagebreak
\medskip 
$g=2$:\nopagebreak
$$\goth h_0 = p^{6}= p^{6}\goth h_4$$
$$\goth h_1 = -p^{3}\tau_p= p^{3}\goth h_3$$
$$\goth h_2 = p(\tau_{p,1} +p^2 + 1) $$
$$\goth h_3 = -\tau_p$$
$$\goth h_4 = 1$$
$g=3$:\nopagebreak
$$\goth h_0 = p^{24}= p^{24}\goth h_8$$
$$\goth h_1 = -p^{18}\tau_p= p^{18}\goth h_7 $$
$$\goth h_2 = p^{13}[\tau_{p,1} +(p^2 + 1)\tau_{p,2} 
+ (-p^5-p^3 + 2p^2 + 1)]= p^{12}\goth h_6$$
$$\goth h_3 = -p^9[\tau_p \tau_{p,2} + \tau_p]= p^{6}\goth h_5$$
$$\goth h_4 = p^6[\tau_p^2 + \tau_{p,2}^2+(-2p+2)\tau_{p,2}-
2p\tau_{p,1}+p^6+2p^4-2p^3-2p+1]$$
$$\goth h_5 = -p^3 [\tau_p \tau_{p,2}+\tau_p]$$
$$\goth h_6 = p[\tau_{p,1} +(p^2 + 1)\tau_{p,2} + (-p^5-p^3 + 2p^2 + 1)]$$
$$\goth h_7 = -\tau_p$$
$$\goth h_8 = 1$$
\medskip
{\bf Table 5. Numbers $\tilde C({\goth d},k)$.}\nopagebreak
\medskip
$\tilde C({\goth d},k)$ depend only on $q_1$ --- the quantity of ones in $\goth d$ --- and $k$; particularly, they don't depend of $g$. 
\medskip
$\matrix & k&0&1&2&3 \\ q_1&&&&&\\0&&1&&&\\ 1&&1-p^{-1}&p^{-1}&&\\ 2&&2-2p^{-1}&p^{-1}-p^{-3}&p^{-3}&\\3&&3-4p^{-1}+p^{-4}&3p^{-1}-p^{-2}-p^{-3}-p^{-4}&p^{-3}-p^{-6}&p^{-6} \endmatrix $
\medskip
{\bf Table 6. Explicit values of $\alpha_G(\chi)(\tau_{p,i})$.}\nopagebreak
\medskip
We denote $\sigma_i=\sigma_i(b_*)$.
\medskip
$g=2$: $\alpha_G(\chi)(\tau_{p,1})= p^{-1}\sigma_1 \sigma_2 + (p^{-1}-p^{-3})\sigma_2 + p^{-1} \sigma_1$
\medskip
$g=3$: $\alpha_G(\chi)(\tau_{p,1})=p^{-1}(\sigma_1+\sigma_1\sigma_2 +\sigma_2\sigma_3) +  
(p^{-1}-p^{-3}) (\sigma_1\sigma_3 + \sigma_2) + (p^{-1}-p^{-4})\sigma_3$
\medskip
$\alpha_G(\chi)(\tau_{p,2})=p^{-3}(\sigma_1\sigma_3 + \sigma_2) + (p^{-3}-p^{-6})\sigma_3$
\medskip
{\bf Table 7. Polynomials $Y_i$.}\nopagebreak
\medskip
$Y_1=-1$

$Y_2=p^3-p^2+p-1$

$Y_3=-p^7+p^6-p^5+p^4+p^3-2p^2+p-1$
\medskip
{\bf 8. Some properties of $\Cal M_P$, $P$ of two simplest types.}\footnote{Written by a request of a reader; can be withdrawn of the text.}\nopagebreak
\medskip
(a) $P=B$, i.e. $P$ of the first type, all $\goth b_i=1$. 
\medskip
In this case $\Cal M_P$ is (generally) irreducible, of weight $\frac{g(g+1)}2$, the packet $\Pi_P$ consists of $2^{g-1}$ representations $\pi_{\goth c}$ where $\goth c$ runs over the set of all subsets of $1,\dots,g$ factorized by the equivalence relation: a subset is equivalent to the complementary subset. The partition (0.1.1) is trivial (i.e. consists of one set). For any $\goth c \in \Pi_P$ $\goth {S(c},P)$ is trivial, and the Hodge number $h^{i,j}(\Cal M_P)$ ($i+j=\frac{g(g+1)}2)$ is equal to the quantity of subsets of $1,\dots,g$ such that the sum of elements of this subset is $i$, i.e. $$\sum_{i=0}^{\frac{g(g+1)}2}h^{i,j}(\Cal M_P)t^i= \prod_{i=1}^g(t^i+1)$$ 
\medskip
(b) $P$ of the second type, all $\goth b_i=1$. 
\medskip
In this case $\Cal M_P$ is the sum of 2 (generally) irreducible submotives $M^-$, $M^+$ of weights $\frac{g(g+1)}2-1$, $\frac{g(g+1)}2+1$ respectively, the packet $\Pi_P$ consists of $2^{g-2}$ representations $\pi_{\goth c}$ where $\goth c$ runs over the set of all subsets of $2,\dots,g$ factorized like in (a). (0.1.1) is also trivial. For any $\goth c \in \Pi_P$ $\goth {S(c},P)$ is the irreducible $\goth{sl}_2$-module of dimension 2, and the Hodge number $h^{i,j}(M^-)$ ($i+j=\frac{g(g+1)}2-1)$ is equal to the quantity of subsets of $2,\dots,g$ such that the sum of elements of this subset is $i$, i.e. $$\sum_{i=0}^{\frac{g(g+1)}2-1}h^{i,j}(M^-)t^i= \prod_{i=2}^g(t^i+1)$$ 
\medskip
{\bf 8a. Hecke polynomial for the case $M^-$, $g=3$.}\nopagebreak
\medskip
We consider numbers $b_1$, $b_2$, $b_3$ for $M^-$. $b_1=p$, we denote $s_1=b_2+b_3$, $s_2=b_2b_3$. Roots of Hecke polynomial for $M^-$ are $a_0$, $a_0b_2$, $a_0b_3$, $a_0b_2b_3$, i.e. this polynomial is
$$\fr^4 - a_0(1+s_1+s_2)\fr^3+ a_0^2 (s_1+s_1s_2+2s_2)\fr^2- a_0 (1+s_1+s_2)p^5 \fr + p^{10}$$
(recall that $a_0^2s_2=p^5$). 

In terms of $\tau_{p}$, $\tau_{p,1}$ this polynomial is 

$\fr^4-\frac{\tau_p}{p+1}\fr^3+\frac{p^2}{(p^2+1)(p-1)} [-(\frac{\tau_p}{p+1})^2 + \tau_{p,1} + p^6 -2 p^5 +2 p^4 -2 p^3 +2 p^2] \fr^2 -\frac{p^5\tau_p}{p+1}\fr + p^{10}$. 

It can be also obtained by taking the Hecke polynomial for $g=3$, substituting $\tau_{p,2}$ from (4.5) and factorizing the obtained expression. 
\medskip
{\bf 9. Another way to find relations between $\goth m_p$, $\goth m_{p,i}$.}\nopagebreak
\medskip
Here we give this method only for the case $P$ of Appendix, 8b. Let $\alpha_1,\dots,\alpha_{2^{g-1}}$ be the Weil numbers of $M^-$, so $p\alpha_1,\dots,p\alpha_{2^{g-1}}$ are the Weil numbers of $M^+$. We denote the Weil numbers of $\Cal M_P$ by $\gamma_1,\dots,\gamma_{2^{g}}$ (the union of sets $\alpha_1,\dots,\alpha_{2^{g-1}}$ and $p\alpha_1,\dots,p\alpha_{2^{g-1}}$), so we have 
$$\sigma_i(\gamma_*)=\sum_{j=0}^i p^j \sigma_{i-j}(\alpha_*) \sigma_{j}(\alpha_*)\eqno{(A1)}$$
Further, $\sigma_i(\gamma_*)=\goth h_i$. From now we consider only the case $g=3$. 

Taking values of $\goth h_i$ from Table 4 and taking into consideration that $\sigma_3(\alpha_*)=p^5\sigma_1(\alpha_*)$, $\sigma_4(\alpha_*)=p^{10}$ we get from (A1), $i=1,2$: 
\medskip
$\sigma_1(\alpha_*)=\frac{\goth m_p}{p+1}$; 

$\sigma_2(\alpha_*)=\frac{p(\goth m_{p,1}+(p^2+1)\goth m_{p,2}-p^5-p^3+2p^2+1-\frac{\goth m_p^2}{(p+1)^2})}{p^2+1}$
\medskip
Further, the equality for $\goth h_3$ is equivalent to $\goth m_p A=0$
and the equality for $\goth h_4$ is equivalent to $AB=0$, where the common multiple $A$ is the left hand side of (4.5): 
\medskip
$A=\frac{\goth m_p^2}{(p+1)^2}+\sum_{i=1}^{3}Y_i\goth m_{p,i}=(p^3-p^2+p-1)\goth m_{p,2}+(-p^7+p^6-p^5+p^4+p^3-2p^2+p-1)+\frac{\goth m_p^2}{(p+1)^2}-
\goth m_{p,1}$ 

and 

$B=-\frac{\goth m_p^2}{(p+1)^2}+\sum_{i=1}^{3}Y_i^+\goth m_{p,i}=(p^3+p^2+p+1)\goth m_{p,2}+(+p^7+p^6+p^5+p^4+p^3+2p^2+p+1)-\frac{\goth m_p^2}{(p+1)^2}+
\goth m_{p,1}$, where $Y_i^+$ are polynomials in $p$ whose coefficients are the absolute values of the ones of $Y_i$. 
\medskip
(4.5) shows that $\goth m_{p,*}$ satisfy the condition $A=0$ (but not $\goth m_{p}=B=0$). 

Substituting the condition $A=0$ to the formula for $\sigma_2(\alpha_*)$, we can slightly simplify it: 
\medskip
$\sigma_2(\alpha_*)=p^2(\goth m_{p,2}-p^4+p^3-p^2+1)$.
\medskip
\medskip
\centerline{{\bf Notation Index}}
\settabs 8 \columns
\+$a_i$&2.1\cr 
\+$\alpha$&After 3.1.2\cr 
\+$\alpha_G(\chi)$&2.1\cr 
\+$\alpha_T(\chi)$&2.5\cr 
\+$\goth a$&3.5\cr 
\+$\goth a_i$&3.5\cr 
\+$b_i$&2.1\cr 
\+$b_I$&0, Step 1\cr 
\+$B$&Borel subgroup, 2.1\cr
\+$B(*)$&Proof of 4.3 \cr 
\+$\goth b_i$&0.2, 3.1.1\cr 
\+$c_i$&3.3\cr 
\+$c_*$&coefficients, after 4.3.6\cr 
\+$\goth c$&3.3\cr 
\+$\goth C$&3.3\cr 
\+$d_i$&3.5\cr 
\+$D$&3.5\cr 
\+$D_i$&3.5\cr 
\+$\goth d$&Line between (4.3.1) and (4.3.2) \cr 
\+$E_g$&the unit matrix\cr
\+$E_{ij}$, $E_{i,j}$&the elementary matrix\cr
\+$\goth e_j$&after 4.3.6\cr 
\+$E$&field of coefficients, Introduction\cr 
\+$f_i$&Remark of 3.5\cr 
\+$F_i$&1.1\cr
\+$\Phi_i$&1.1\cr
\+$g$&genus\cr 
\+$G$&$GSp_{2g}$\cr 
\+$\gamma$&3.2\cr
\+$h$&0; 3.4\cr 
\+$\n H(*)$&1.1\cr 
\+$\goth h_i$&1.5\cr 
\+$\theta_{\pi_{\chi}}$&Langlands element 2.7\cr 
\+$k$&0.2, 3.1.1\cr 
\+$K$&3.5\cr
\+$K_c$&Centralizer of $\mu$ 3.5\cr 
\+$l$&length of an element of $\Omega$ \cr 
\+$m_D$&3.5\cr 
\+$m_{D_i}$&3.5\cr 
\+$m_i$, $m_j$&after 0.2; after 3.1.1\cr 
\+$m_{\rho}$&3.5\cr 
\+$M_s$&1.1\cr
\+$M$&Levi subgroup 3.1\cr 
\+$M_c$&Levi subgroup after 3.1.2\cr 
\+$\mu$&3.4\cr 
\+$\Cal M$&Introduction\cr 
\+$\Cal M_P$&after 0.1\cr 
\+$\goth m$&after 0.1\cr
\+$\goth m_p$, $\goth m_{p,i}$&above 4.4\cr
\+$\goth m_{\Cal M}$&Introduction\cr
\+$N$&3.1\cr
\+$N_c$&after 3.1.2\cr
\+$\Cal N_*$&3.3\cr
\+$\goth N$&normalizer\cr
\+$\Omega$&Weyl group\cr 
\+$p$&prime, Introduction\cr 
\+$p_w$&3.4\cr 
\+$\pi_{\goth c}$, $\pi_w$&3.3\cr 
\+$\pi_{\chi}$&2.1\cr 
\+$P$&parabolic subgroup 3.1\cr 
\+$\Pi_P$&0, p. 3; 3.3\cr 
\+$q_w$&3.4\cr 
\+$r$&2.6\cr 
\+$R_n(i)$&1.3\cr
\+$\rho$&3.5\cr 
\+$S(n)$&permutations group\cr 
\+$S_*$&Satake map, 1.1\cr
\+$\sigma_i$&$i$-th symmetric polynomial\cr 
\+$\goth S(c,b)$&0, Step 2c; 3.5\cr 
\+$\goth S(\goth c,P)$&0.3; 3.5\cr 
\+$T$&diagonal in $G$, 1.1\cr
\+$T_{p}$, $T_{p,i}$&1.1\cr
\+$\tau_{p}$, $\tau_{p,i}$&generators of Hecke algebra, 1.1\cr 
\+$U_i$&1.5\cr
\+$U_I$&1.5\cr
\+$V_i$&1.5\cr
\+$V_{\pi}$&after 4.3\cr 
\+$V_{P}$&after 4.3\cr 
\+$w$&3.3.1, 3.3.2\cr 
\+${\goth w}$&3.4\cr 
\+${\goth w}_I$&end of 3.2\cr 
\+$W$&2.6\cr 
\+${\goth W}$&3.3\cr 
\+$x_I$&2.6\cr 
\+$x_n$&Between 4.3.6 and 4.3.7\cr 
\+$X$&Shimura variety, p. 2 \cr
\+$\goth x_{\goth c}$&4.1\cr
\+$\chi$&2.1\cr
\+$Y_i$&4.6\cr
\medskip
\medskip
{\bf References}\nopagebreak
\medskip
[AZh] Andrianov A.N., Zhuravlev V.G., Modular forms and Hecke
operators. Moscow, 1990 (In Russian; English version: 
Andrianov A.N., Quadratic forms and Hecke operators. Springer, 1987).

[A] Arthur J., Unipotent automorphic representations: conjectures. 
Asterisque, 171 - 172 (1989), p. 13 - 71

[AJ] Adams J., Johnson J.F., Endoscopic groups and packets of non-
tempered 
representations. Compositio Math., 64, 271 -309 (1987)

[BR] Blasius D., Rogawski J.D., Zeta functions of Shimura varieties.
In: Motives. Proc. of Symp. in Pure Math., 1994, v. 55, part 2, p. 
525 - 571. 

[B] Bultel O., On the mod $\goth P$-reduction of 
ordinary CM-points. Oxford, 1997. Ph. D. thesis.  

[D] Deligne P. Travaux de Shimura. Lect. Notes in Math., 1971, 
v.244, p. 123 - 165. Seminaire Bourbaki 1970/71, Expos\'e 389. 

[FCh] Faltings G., Chai Ching-Li. Degeneration of 
abelian varieties. Springer, 1990

[J] Jacobson. Lie algebras. Interscience tracts in pure and applied mathematics, vol. 10. 1962 

[S] Satake Ichiro, Theory of spherical functions 
on reductive algebraic groups over $\goth p$-adic fields. 
Publ. IHES, 1963. p. 229 - 293

[VZ] Vogan D., Zuckerman G. Unitary representations with non-zero cohomology. Compositio Math., 53 (1984), 51 -- 90 
\medskip
\medskip
{\bf Address. }
\medskip
Departamento de Matem\'aticas 

Universidad Sim\'on Bol\'\i var

Apartado Postal 89000

Caracas, Venezuela
\medskip
Tel. (58)(212)9063283 

Fax (58)(212)9063278 
\medskip
E-mail: logachev\@usb.ve 
\end